\documentclass[final,leqno,onefignum,onetabnum]{siamltex1213}
\usepackage{amsfonts}
\usepackage{amsmath,booktabs,ctable,threeparttable}
\usepackage{amssymb,amsfonts,boxedminipage}
\usepackage{algorithm}
\usepackage{algorithmic}

\title{Modified Truncated Randomized Singular Value Decomposition (MTRSVD)
Algorithms for Large Scale Discrete Ill-posed Problems with General-Form
Regularization\thanks{This
work was supported in part by
the National Science Foundation of China (Nos. 11771249 and 11371219)}}

\author{Zhongxiao Jia\thanks{Corresponding author.
Department of Mathematical Sciences, Tsinghua
University, 100084 Beijing, China. (\email{jiazx@tsinghua.edu.cn})} \and
Yanfei Yang\thanks{Department of Mathematical Sciences, Tsinghua
University, 100084 Beijing, China. (\email{yangyf14@mails.tsinghua.edu.cn})}}

\begin{document}
\maketitle
\slugger{sirev}{xxxx}{xx}{x}{x--x}

\begin{abstract}
In this paper, we propose new randomization based algorithms for large scale linear
discrete ill-posed problems with general-form regularization:
${\min} \|Lx\|$ subject to ${\min} \|Ax - b\|$, where $L$ is a regularization matrix.
Our algorithms are inspired by the modified truncated singular
value decomposition (MTSVD) method, which suits only for small to medium scale problems,
and randomized SVD (RSVD) algorithms that generate good low rank approximations to $A$.
We use rank-$k$ truncated randomized SVD (TRSVD) approximations to $A$
by truncating the rank-$(k+q)$ RSVD approximations to $A$, where $q$ is an
oversampling parameter. The resulting algorithms are called
modified TRSVD (MTRSVD) methods. At every step, we use the LSQR
algorithm to solve the resulting inner least squares problem, which is proved to
become better conditioned as $k$ increases so that LSQR converges faster.
We present sharp bounds for the approximation accuracy of the RSVDs and TRSVDs
for severely, moderately and mildly ill-posed problems, and
substantially improve a known basic bound for TRSVD approximations. We prove
how to choose the stopping tolerance for LSQR in order to
guarantee that the computed and exact best regularized solutions
have the same accuracy. Numerical experiments illustrate
that the best regularized solutions by MTRSVD are as accurate as
the ones by the truncated generalized singular
value decomposition (TGSVD) algorithm, and at least as accurate
as those by some existing truncated
randomized generalized singular value decomposition (TRGSVD) algorithms.
\end{abstract}



\pagestyle{myheadings}
\thispagestyle{plain}
\markboth{Z. JIA AND Y. YANG}{MTRSVD ALGORITHMS FOR ILL-POSED PROBLEMS}

\section{Introduction}

Consider the solution of the large-scale linear discrete ill-posed problem
\begin{equation}\label{eq1}
\min_{x\in \mathbb{R}^{n}}\|Ax - b\| \quad {\rm or}
\quad Ax = b, \quad A \in \mathbb{R}^{m \times n}, \quad b \in \mathbb{R}^{m},
\end{equation}
where the norm $\|\cdot\|$ is the 2-norm of a vector or matrix,
the matrix $A$ is ill conditioned with its singular values
decaying to zero with no obvious gap between consecutive ones,
and the right-hand side $b = b_{true}+e$ is noisy and
assumed to be contaminated by a white noise $e$, which may stem from
measurement, truncation or discretization errors, where $b_{true}$
represents the unknown noise-free right-hand side and
$\|e\| < \|b_{true}\|$.
Such kind of problem arises in a variety of applications, such as
computerized tomography, electrocardiography, image deblurring,
signal processing, geophysics, heat
propagation, biomedical and optical imaging, groundwater modeling,
and many others; see, e.g.,
\cite{aster,cuppen,elden82,engl93,kaipio,miller,nat1986}.

The naive solution $x_{naive} = A^{\dag}b$
is a meaningless approximation to the true solution
$x_{true} = A^{\dag}b_{true}$ since $b$ is
contaminated by the noise and $A$ is extremely ill conditioned, where $\dag$
denotes the Moore-Penrose inverse of a matrix. Therefore, one has
to use regularization to obtain a best possible approximation to $x_{true}$
\cite{hansen98,hansen10}.

One of the common regularization approaches is to solve the standard-form
regularization problem
\begin{equation}\label{common}
\min\|x\|\quad{\rm subject \ \ to}
\quad\|Ax-b\|=\min.
\end{equation}
The truncated singular value decomposition (TSVD) method is one of
the most popular regularization methods for solving \eqref{common}.
The method computes a minimum 2-norm
least squares solution, i.e., the TSVD solution $x_k$,  which
solves the problem
\begin{equation}\label{tik}
\min\|x\|\quad{\rm subject \ \ to}
\quad x\in \mathcal{S}_k = \left\{x\mid\|A_kx-b\|=\min\right\}
\end{equation}
starting with $k=1$ onwards until a best regularized solution is found
at some $k$, where $A_k$ is a best rank-$k$ approximation to $A$
with respect to the 2-norm and the index $k$ plays the role of the regularization
parameter. It is known from, e.g., \cite[p.~79]{golub13}, that
\begin{align}\label{best}
\|A-A_k\|=\sigma_{k+1},
\end{align}
where $\sigma_{k+1}$ is the $(k+1)$th large singular value of $A$.
\eqref{common} is equivalent to the standard-form Tikhonov regularization problem
\begin{align}\label{general2}
\min_{x\in \mathbb{R}^n}\left\{\|Ax-b\|^2 +\lambda^2\|x\|^2\right\}
\end{align}
with the regularization parameter $\lambda>0$.  \eqref{common} amounts to
\eqref{general2} in the sense that for any regularization parameter
$\lambda\in[\sigma_n,\sigma_1]$ there is a truncation
parameter $k$ such that the solutions computed by the TSVD method and
the Tikhonov regularization method  are close.
Furthermore, with the optimal parameter $\lambda_{opt}$ chosen,
the best regularized solutions obtained by the
two methods have very comparable
accuracy with essentially the minimum 2-norm error \cite{hansen98,hansen10}.

Hansen \cite{hansen98} points out that, in many applications,
minimizing the $2$-norm of the solution, i.e., $\min\|x\|=\min\|I_nx\|$
with $I_n$ being the $n\times n$ identity matrix, is not an optimal choice.
On the one hand, $\|x\|$ may not always be affected as much by the
errors as the $2$-norm of a derivative of the solution.
On the other hand, the SVD basis vectors
may not be well suited for computing a good regularized solution
to \eqref{eq1}, but choosing a regularization matrix $L\neq I_n$
can often lead to a much better approximate solution.
He presents some examples such as data approximation by bivariate spline \cite{cox}.
Kilmer {\em et al}. \cite{hansen2007} also give some examples from
geophysics and heat distribution, where choosing an $L\neq I_n$ appears
more effective.

In this paper, we consider to exploit
the priori information on $x_{true}$ by using $\min\|Lx\|$ in \eqref{common}
and \eqref{tik} other than $\min\|x\|$, that is, we solve the
general-form regularization problem
\begin{equation}\label{general}
\min\|Lx\| \mbox{ \ \
subject to\ } \|Ax-b\|=\min,
\end{equation}
where $L\in\mathbb{R}^{p\times n}$
is usually a discrete approximation of some derivative operators.
When $L \neq I_n$, \eqref{general2} becomes the general-form Tikhonov
regularization problem
\begin{align}\label{general3}
\min_{x\in\mathbb{R}^{n}}\left\{\|Ax-b\|^2+\lambda^2\|Lx\|^2\right\},
\end{align}
which is equivalent to \eqref{general}.
The solution to \eqref{general3} is unique for a given $\lambda>0$ when
$$\mathcal{N}(A)\cap \mathcal{N}(L)={0}\Longleftrightarrow
\rank\left(
\begin{array}{c}
A \\
L \\
\end{array}
 \right)=n,$$
where $\mathcal{N}(\cdot)$ denotes the null space of a matrix.
In practical applications, $L$ is typically chosen as
\begin{equation}\label{l1}
L_1 = \left(
        \begin{array}{ccccc}
          1 & -1 &  &  &  \\
           & 1 & -1 &  &  \\
           &  & \ddots & \ddots &  \\
             &  &  & 1  & -1\\
        \end{array}
      \right)\in \mathbb{R}^{(n-1)\times n},
\end{equation}

\begin{equation}\label{l2}
L_2 = \left(
        \begin{array}{cccccc}
          -1 & 2 & -1 &  &  &\\
           & -1 & 2 & -1 &  &\\
           &  & \ddots & \ddots & \ddots & \\
           &  &  & -1 &2  &-1\\
        \end{array}
      \right)\in \mathbb{R}^{(n-2)\times n},
\end{equation}
or
\begin{equation}\label{l1l2}
L_3 = \left(
        \begin{array}{c}
          L_1 \\
          L_2 \\
        \end{array}
      \right)\in \mathbb{R}^{(2n-3)\times n},
\end{equation}
where $L_1$ and $L_2$ are scaled discrete approximations of
the first and second derivative operators in one dimensional
Fredholm integral equations of the first kind, respectively.
For the corresponding
$L_1$ and $L_2$ in two dimensional problems, see Section 8.2
of \cite{hansen10}.

For small to medimum scale problems, adapting the TSVD method
to Problem \eqref{general},
Hansen {\em et al.} \cite{hansen92} propose a
modified truncated SVD (MTSVD) method that solves
\begin{equation}\label{mtik}
\min\|Lx\|\quad {\rm subject \ \ to}\quad x\in\mathcal{S}_k
= \left\{x\mid\|A_kx-b\|=\min\right\}
\end{equation}
starting with $k=1$ onwards until a best regularized solution is found for
some $k$. As in the TSVD method, $k$ plays the role of
the regularization parameter. This approach is an alternative to the
TGSVD method for solving \eqref{general3}.
The algorithm first computes the SVD of $A$ and then extracts the
best rank-$k$ approximation $A_k$ to $A$ by truncating the SVD of $A$.
It solves a sequence of least squares problems by the adaptive
QR factorization from $k=1$ onwards until
a best regularized solution is found. This algorithm avoids
computing the GSVD of the matrix pair $\{A, L\}$, but it is not
suitable for large scale problems since computing
the SVD of $A$ is infeasible for $A$ large.

For $L=I_n$, Xiang and Zou \cite{zou13} adapt some basic
randomized algorithms from \cite{Rev2011} to \eqref{general2} and
develop a randomized SVD (RSVD) algorithm. RSVD
acts $A$ on a Gaussian random matrix to capture the
dominant information on the range of $A$, and
computes the SVD of a small matrix. By the SVD of the small matrix, one
then obtains an approximate SVD of $A$. Halko {\em et al}. \cite{Rev2011} have
given an accuracy analysis on the randomized algorithm and
approximate SVD, and have established a number of error bounds for them.
Randomized algorithms have been receiving high attention in recent years and
widely used in a variety of low rank approximations; see, e.g.,
\cite{gu15, gu16, Rev2011, liberty, mar, lsrn, PCA, wei16, zou13, zou15}.

For $L\neq I_n$, Xiang and Zou \cite{zou15} present a randomized
GSVD (RGSVD) algorithm to solve \eqref{general3}.
First, they compute a RSVD of $A$. Then they compute the GSVD of
the matrix pair $\{AQ,LQ\}$, where $Q$
is the right singular vectors in RSVD. The matrix $Q$
captures the information on dominant right singular vectors of $A$,
which ensures that $AQ$ captures the dominant left singular vectors
of $A$. Indeed, $AQQ^T$ is a good low rank approximation to $A$
with high probability; see \cite{Rev2011} for some bounds and next section
for refined bounds. However, the generation
of $Q$ does not make use of any information on $L$. As a consequence, there is no
guarantee that the GSVD of
the matrix pair $\{AQ,LQ\}$ well approximates the dominant GSVD components of
$\{A,L\}$, which is a critical requirement that RGSVD can obtain a good
regularized solution to \eqref{eq1}.

Wei {\em et al}. \cite{wei16} propose new RGSVD algorithms.
For the underdetermined case, their algorithm is the
same as Xiang and Zou \cite{zou15} in theory. An algorithmic difference is that
they do not compute an approximate SVD of $A$. Instead,
they compute the GSVD of the matrix pair $\{AQ,LQ\}$,
where $Q$ captures only the information on dominant right singular vectors
of $A$ and has nothing to do with $L$. Therefore, it has the same deficiency
as the algorithm in \cite{zou15}, as mentioned above. For the
overdetermined case, their RGSVD method needs to compute the GSVD of
the matrix pair $\{B,L\}$, where $B=Q^TA\in\mathbb{R}^{l\times n}$ is
a dense matrix with $Q$ being an $m\times l$ orthonormal matrix
generated by randomized algorithms, $L\in\mathbb{R}^{p\times n}$
and the parameter $l$ satisfies $l+p\geq n$.
Since this algorithm captures the dominant information on $A$ and
retains $L$ itself, it works theoretically for \eqref{eq1}.
However, for a large scale \eqref{eq1}, $n$ must be large, so is
the size of the matrix pair $\{B,L\}$. This leads to
the computation and storage memory of the GSVD of the matrix
pair $\{B,L\}$ impractical because one must
compute a large dense $n\times n$ matrix and invert it to obtain the
right singular vector matrix of this matrix pair when
using the resulting RGSVD to solve \eqref{general} or \eqref{general3}.
As a result, the proposed RGSVD algorithm actually does not suit
for large scale problems.

In this paper, inspired by the idea of randomized algorithms and the MTSVD method,
we will propose a modified truncated randomized SVD (MTRSVD) method for
solving \eqref{general}.
Our method consists of four steps: first, use the RSVD algorithms \cite{Rev2011}
to obtain approximate SVDs of $A$ for the underdetermined and overdetermined cases,
respectively; second, truncate the approximate SVDs to
obtain rank-$k$ TRSVD approximations $\widetilde{A}_k$ to $A$;
third, use $\widetilde{A}_k$ to replace the best rank-$k$
approximation $A_k$ in \eqref{mtik};
finally, solve
\begin{equation}\label{mtik2}
\min\|Lx\| \quad {\rm subject\ \ to}\quad x\in\mathcal{S}_k
=\{x\mid\|\widetilde{A}_kx-b\|=\min \}
\end{equation}
starting with $k=1$ onwards until a best regularized solution is found for
some $k$. As will be seen later, this step gives rise to a large least squares
problem that is different from the one in \cite{hansen92} and cannot
be solved by adaptive QR factorizations any more because of its large
size and the unavailability of the SVD of $A$.
We will use the LSQR algorithm \cite{paige1982}
to iteratively solve the resulting least squares problem.

We consider a number of theoretical issues on the MTRSVD algorithms.
For severely, moderately and mildly ill-posed problems \cite{hansen98,hansen10,hof},
we establish some sharp error bounds for the approximation error $\|A-QQ^TA\|$
(or $\|A-AQQ^T\|$)
in terms of $\sigma_{k+1}$, where $QQ^TA$ (or $AQQ^T$) is
the rank-$(k+q)$ RSVD approximation and $Q\in \mathbb{R}^{m\times (k+q)}$
(or $Q\in \mathbb{R}^{n\times (k+q)})$ is an
orthonormal matrix with $q$ being an oversampling parameter.
Halko {\em et al.} \cite{Rev2011} have presented a number of error bounds
for the approximation errors.
Their bounds have been used in, e.g., \cite{wei16,zou15} and are good enough for
a nearly rank deficient $A$, but turn out to be possibly meaningless for
ill-posed problems since they are pessimistic and even may never become small
for any $k$ and $q$. In contrast, our bounds are always meaningful
and much sharper for the aforementioned three kinds of ill-posed problems.
Next, for the truncated rank-$k$ approximations
$\widetilde{A}_k$, we focus on a basic bound in \cite{Rev2011} and improve
it substantially. Our new bounds are unconditionally superior to
and can be much sharper than the bound for $\widetilde{A}_k$ in \cite{Rev2011},
and they explain why the error introduced in truncation step is not
so damaging, an important concern in \cite[Remark 9.1]{Rev2011}. For the
MTRSVD algorithms, we analyze
the conditioning of the resulting inner least squares problem at each step $k$.
We will prove that the condition number monotonically decreases as $k$ increases,
such that  for the same stopping tolerance the LSQR algorithm for solving
it generally converges faster
and uses fewer inner iterations as $k$ increases.
In the meantime, we consider efficient implementations of
Lanczos bidiagonalization used within LSQR for the inner least squares problems.
Importantly, we will make a detailed analysis on the stopping tolerance for LSQR,
showing how to choose it so as to guarantee that the computed and exact best
regularized solution have the same accuracy.
We prove that the stopping tolerance
for LSQR is not restrictive and a reasonably small one is good enough,
provided that the regularization matrix $L$ is well conditioned.
Finally, we report numerical experiments to illustrate the generality
and effectiveness of our algorithms.
We show that, for the $m\geq n$ case
with $n$ not large, the best regularized solutions obtained by MTRSVD
are as accurate as those by the TGSVD algorithm and the TRGSVD algorithm
in \cite{wei16}. When $n$ is large, the TRGSVD algorithm in \cite{wei16}
is out of memory in our computer, but
our algorithm works well. For the $m\leq n$ case, the best
regularized solutions by MTRSVD are very comparable to those by
the TGSVD algorithm and are at least as accurate as those by the
TRGSVD algorithms in \cite{wei16,zou15}.

Our paper is organized as follows. In Section 2, we review the RSVD algorithms
and establish new error bounds for the RSVD approximations to $A$ for
severely, moderately and mildly ill-posed problems, respectively.
In Section 3, we present the MTRSVD algorithms, establish new sharp bounds for
the TRSVD approximation to $A$, and make an analysis on the conditioning
of inner least squares problems and on the stopping tolerance for LSQR.
In Section 4, we report numerical examples to
illustrate that our algorithms work well. Finally, we conclude the paper in Section 5.

\section{RSVD and sharp error bounds}

Let the compact SVD of $A \in\mathbb{R}^{m\times n}$ be
\begin{align}\label{svd}
A = U\Sigma V^T,
\end{align}
where
$U=\left(u_1, u_2,\ldots, u_s\right)\in\mathbb{R}^{m\times s}$ and $
V=\left(v_1, v_2,\ldots, v_s\right)\in\mathbb{R}^{n\times s}$
are column orthonormal, $\Sigma = \diag(\sigma_1, \sigma_2, \ldots, \sigma_s)
\in\mathbb{R}^{s\times s}$ with $s = \min\{m,n\}$ and
$\sigma_1, \sigma_2, \ldots, \sigma_s$ being the
singular values and labeled as $
\sigma_1\geq \sigma_2 \geq \cdots \geq \sigma_s>0
$. Then
\begin{align}\label{ak}
A_k = U_k\Sigma_kV_k^T
\end{align}
is one of the best rank-$k$ approximations to $A$ with respect to the 2-norm,
where $U_k=\left(u_1, u_2,\ldots, u_k\right)\in\mathbb{R}^{m\times k}$ and $V_k=
\left(v_1, v_2,\ldots, v_k\right)\in\mathbb{R}^{n\times k}$ are column orthonormal,
and $\Sigma_k=\diag(\sigma_1, \sigma_2, \ldots, \sigma_k)\in\mathbb{R}^{k\times k}$.
Define the condition number of $A$ as
$$
\kappa(A) = \frac{\sigma_{\max}(A)}{\sigma_{\min}(A)}=\frac{\sigma_1}{\sigma_s}.
$$

The following Algorithm~\ref{alg1} is the basic randomized algorithm,
presented in \cite{Rev2011}, that computes a low rank approximation
to $A$ and an approximate
SVD of $A$ for the overdetermined case ($m\geq n$).

\begin{algorithm}[htb]
\caption{(RSVD) Given $A\in \mathbb{R}^{m\times n}(m\geq n)$, $l=k+q<n$ and $q\geq4$,
compute an approximate SVD: $A\approx \widetilde{U}\widetilde{\Sigma }\widetilde{V}^T$
with $\widetilde{U}\in \mathbb{R}^{m\times l}$, $\widetilde{\Sigma }\in
\mathbb{R}^{l\times l}$ and $\widetilde{V}\in\mathbb{R}^{n\times l}$.}
\begin{algorithmic}[1]\label{alg1}
\STATE Generate an $n\times l$ Gaussian random matrix $\Omega$.
\STATE Form the $m\times l$ matrix $Y=A\Omega$.
\STATE Compute the $m\times l$ orthonormal matrix $Q$ via QR factorization $Y = QR$.
\STATE Form the $l\times n$ matrix $B = Q^TA$.
\STATE Compute the compact SVD of the small matrix $B$:
$B=W\widetilde{\Sigma}\widetilde{V}^T$.
\STATE Form the $m\times l$ matrix $\widetilde{U} =QW$, and
$A\approx \widetilde{U}\widetilde{\Sigma}\widetilde{V}^T= QQ^TA$.
\end{algorithmic}
\end{algorithm}

The mechanism of Algorithm \ref{alg1} is as follows:
the information of the column space of $A$ is extracted in step 2 , i.e.,
$\mathcal{R}(Y)\subseteq \mathcal{R}(A)$ where $\mathcal{R}(\cdot)$
denotes the column space or range of a matrix. It is clear that the
columns of $Q$ span the main range of $A$ in step 3 and
$\mathcal{R}(Q)=\mathcal{R}(Y)\subseteq\mathcal{R}(A)$.
Noting $\mathcal{R}(A)=\mathcal{R}(U)$, the factor $Q$ captures the
dominant left singular vectors of $A$. In step 4, because of
$\mathcal{R}(B^T)\subseteq\mathcal{R}(A^T)=\mathcal{R}(V)$,
the matrix $B$ provides information on the dominant
right singular vectors of $A$. In step 6, the algorithm modifies the approximate
left singular vectors.

For the underdetermined case ($m\leq n$), Halko {\em et al.} \cite{Rev2011}
present Algorithm \ref{alg2}, which
is equivalent to applying Algorithm~\ref{alg1} to $A^T$.

\begin{algorithm}[htb]
\caption{(RSVD) Given $A\in \mathbb{R}^{m\times n}(m\leq n)$, $l=k+q<m$
and $q\geq4$, compute an approximate SVD: $A\approx \widetilde{U}
\widetilde{\Sigma}\widetilde{V}^T$ with $\widetilde{U}\in \mathbb{R}^{m\times l}$,
$\widetilde{\Sigma}\in \mathbb{R}^{l\times l}$ and $\widetilde{V}\in\mathbb{R}^{n\times l}$.}
\begin{algorithmic}[1]\label{alg2}
\STATE Generate an $l\times m$ Gaussian random matrix $\Omega$.
\STATE Form the $l\times n$ matrix $Y=\Omega A$.
\STATE Compute the $n\times l$ orthonormal matrix $Q$ via QR factorization $Y^T = QR$.
\STATE Form the $m\times l$ matrix $B = AQ$.
\STATE Compute the compact SVD of the small matrix $B$:
$B=\widetilde{U}\widetilde{\Sigma}W^T$.
\STATE Form the $n\times l$ matrix $\widetilde{V} =QW$, and
$A\approx \widetilde{U}\widetilde{\Sigma}\widetilde{V}^T=AQQ^T$.
\end{algorithmic}
\end{algorithm}

When $q\geq 4$, Halko {\em et al}. \cite{Rev2011} establish the following basic
estimate on the approximation accuracy of $QQ^TA$ generated by
Algorithm \ref{alg1}:
\begin{align}\label{error1}
\|A - QQ^TA\| \leq \left(1+6\sqrt{(k+q)q{\rm log}q}\right)\sigma_{k+1}+3\sqrt{k+q}
\left(\sum_{j>k}\sigma_j^2\right)^{1/2}
\end{align}
with failure probability at most $3q^{-q}$.
Based on \eqref{error1}, Halko {\em et al}. \cite{Rev2011} derive
a simplified elegant error bound
\begin{equation}\label{error2}
\|A - QQ^TA\| \leq \left(1+9\sqrt{(k+q)(n-k)}\right)\sigma_{k+1}
\end{equation}
with failure probability at most $3q^{-q}$. As we can see clearly, for a fixed $k$
the above two bounds monotonically increases with $q$, which is not in accordance
with a basic result that, for a fixed $k$, the left hand side of \eqref{error2}
monotonically decreases with $q$; see Proposition 8.5 of \cite{Rev2011}.
Xiang and Zou \cite{zou13} and Wei {\em et al}. \cite{wei16}
directly exploit the bound \eqref{error2} in their analysis.
For a nearly rank deficient $A$,
the monotonic increasing property of the right-hand sides of \eqref{error1} and
\eqref{error2} with $q$
do not have serious harm since the bound can be small enough to detect
the numerical rank $k$ whenever $q$ is not large, the singular values
$\sigma_k\gg \sigma_{k+1}$ and $\sigma_{k+1}$ is numerically small.
In the context of ill-posed problems,
however, the situation is completely different since the
bound \eqref{error2} may be too pessimistic and meaningless, as will be clear soon.

We notice another basic bound from \cite{Rev2011}
that has received little attention but appears more insightful and useful
than \eqref{error1}, at least in the context of ill-posed problems:
\begin{equation}\label{error3}
\|A - QQ^TA\| \leq \left(1+16\sqrt{1+\frac{k}{q+1}}\right)\sigma_{k+1}+
\frac{8\sqrt{k+q}}{q+1}\left(\sum_{j>k}\sigma_j^2\right)^{\frac{1}{2}}.
\end{equation}
with failure probability at most $3e^{-q}$.

In the manner of deriving \eqref{error2} from \eqref{error1}, we have a simplified form
of \eqref{error3}:
\begin{equation}\label{overest}
\|A - QQ^TA\| \leq \left(1+16\sqrt{1+\frac{k}{q+1}}+
\frac{8\sqrt{(k+q)(n-k)}}{q+1}\right)\sigma_{k+1}.
\end{equation}
On contrary to \eqref{error1}
and \eqref{error2}, an advantage of the bounds \eqref{error3} and \eqref{overest}
is that they monotonically decrease with $q$ for a given $k$.
Compared with \eqref{error1},
a minor theoretical disadvantage of \eqref{error3} is that its failure probability
$3e^{-q}$ is a little higher than $3q^{-q}$ of \eqref{error1} for $q\geq 4$.
But this should not cause any essential problem for practical purposes.

For the bound \eqref{overest},
it is easily justified that the factor in front of
$\sigma_{k+1}$ lies between $1+\mathcal{O}(\sqrt{n})$ and $1+\mathcal{O}(n)$
for a small fixed $q$, and it is $1+\mathcal{O}(\sqrt{n})$ for $k$ not big when
dynamically choosing $q=k$ roughly. In contrast, for the bound \eqref{error2},
the factor in front of $\sigma_{k+1}$ ranges from
$1+\mathcal{O}(\sqrt{n})$ to $1+\mathcal{O}(n)$ for any $q$.

We will show that the bounds \eqref{error2} and \eqref{overest} may be
fatal overestimates in the context of ill-posed problems.
Based on \eqref{error3} and following Jia's works \cite{Jia17,Jia17(2)},
we carefully analyze the approximation accuracy of $QQ^TA$
for three kinds of ill-posed problems:
severely, moderately and mildly
ill-posed problems, and establish much more accurate  bounds.

Before proceeding, we first give a precise characterization of the degree of
ill-posedness of \eqref{eq1} which was introduced in \cite{hof} and has been widely
used in, e.g., the books \cite{aster, engl2000, hansen98, hansen10, mueller2012}.

\begin{definition}
 If
 $\sigma_j=\mathcal{O}(\rho^{-j}),\ j=1,2,\ldots,n
 $
 with $\rho>1$, then \eqref{eq1} is severely ill-posed.
 If the singular values $
 \sigma_j=\mathcal{O}(j^{-\alpha}), \ j=1,2,\ldots,n,
 $
 then \eqref{eq1} is mildly or moderately ill-posed for $\frac{1}{2}<\alpha\leq1$
 or $\alpha>1$.
\end{definition}

We mention that the requirement $\alpha>\frac{1}{2}$ does not appear in
the aforementioned books but it is added in \cite{Jia17,Jia17(2)},
where it is pointed out that this requirement is naturally met
when the kernel of an underlying linear
Fredholm equation of the first kind is square integrable over
a defined domain.

Keep in mind that the factors in front of $\sigma_{k+1}$ in
\eqref{error2} and \eqref{overest} lie
between $1+\mathcal{O}(\sqrt{n})$ and $1+\mathcal{O}(n)$ for a given
$k$. However, for moderately and mildly ill-posed problems,
the bounds \eqref{error2} and \eqref{overest} may never be small
for $k$ not big and $\alpha$ close to one; for $\alpha$ close to
$\frac{1}{2}$, they are definitely not small as $k$ increases up
to $n-q-1$. These bounds, if realistic,
mean that Algorithm \ref{alg1} may never generate
a meaningful rank-$(k+q)$ approximation to $A$. Fortunately, as we will
show below, the bound \eqref{overest} can be improved substantially for
the three kinds of ill-posed problems, and the new bounds indicate that
Algorithm~\ref{alg1} (or Algorithm~\ref{alg2}) indeed generates
very accurate rank-$(k+q)$ approximations to $A$.

\begin{theorem}\label{th1}
For the severely ill-posed problems with
$\sigma_j=\mathcal{O}(\rho^{-j})$ and $\rho>1$, $ j=1,2,\ldots,n,$
it holds that
\begin{align}\label{error5}
\|A - QQ^TA\| \leq \left(1+16\sqrt{1+\frac{k}{q+1}}+\frac{8\sqrt{k+q}}{q+1}
\left(1+\mathcal{O}(\rho^{-2})\right)\right)\sigma_{k+1}
\end{align}
with failure probability at most $3e^{-q}$ for $q\geq 4$ and $k=1,2,\ldots, n-q-1$.
\end{theorem}

{\em Proof}.
By the assumption on the singular values $\sigma_j$, we obtain
\begin{align}
\left(\sum_{j=k+1}^n\sigma_j^2\right)^{1/2}
&= \sigma_{k+1} \left(\sum_{j=k+1}^n
\frac{\sigma_j^2}{\sigma_{k+1}^2}\right)^{1/2}
=\sigma_{k+1} \left(1+\sum_{j=k+2}^n
\frac{\sigma_j^2}{\sigma_{k+1}^2}\right)^{1/2} \notag\\
&=\sigma_{k+1}\left(1+\sum_{j=k+2}^n\mathcal{O}
(\rho^{2(k-j)+2})\right)^{1/2}
\notag \\
&=\sigma_{k+1}\left(1+\mathcal{O}\left(\sum_{j=k+2}^n
\rho^{2(k-j)+2}\right)\right)^{1/2}
\notag \\
&=\sigma_{k+1}\left(1+ \mathcal{O}\left(\frac{\rho^{-2}}
    {1-\rho^{-2}}\left(1-\rho^{-2(n-k-1)}\right)\right)\right)^{1/2}\notag \\
&=\sigma_{k+1} \left(1+\mathcal{O}(\rho^{-2})\right)^{1/2}\notag\\
&=\sigma_{k+1} \left(1+\mathcal{O}(\rho^{-2})\right).
\label{severe1}
\end{align}
Substituting \eqref{severe1} into \eqref{error3} gives \eqref{error5}.
\qquad\endproof

\begin{theorem}\label{th2}
For the moderately and mildly ill-posed problems with
$\sigma_j=\zeta j^{-\alpha}$, $j=1,2,\ldots,n$, where
$\alpha>1/2$ and $\zeta>0$ is some
constant, it holds that
\begin{align}\label{error4}
\|A - QQ^TA\| \leq \left(1+16\sqrt{1+\frac{k}{q+1}}+
\frac{8\sqrt{k+q}}{q+1}\sqrt{\frac{k}{2\alpha-1}}\left(\frac{k+1}{k}\right)^{\alpha}
\right)\sigma_{k+1}
\end{align}
with failure probability at most $3e^{-q}$
for $q\geq 4$ and $k=1,2,\ldots, n-q-1$.
\end{theorem}

{\em Proof}.
By the assumption on the singular values $\sigma_j$, we obtain
\begin{align}
\left(\sum_{j=k+1}^n\sigma_j^2\right)^{1/2}
&= \sigma_{k+1} \left(\sum_{j=k+1}^n
\frac{\sigma_j^2}{\sigma_{k+1}^2}\right)^{1/2}\notag\\
&= \sigma_{k+1} \left(\sum_{j=k+1}^n \left(\frac{j}{k+1}
\right)^{-2\alpha}\right)^{1/2} \notag \\
&=\sigma_{k+1}
\left((k+1)^{2\alpha}\sum_{j=k+1}^n \frac{1}{j^{2\alpha}}\right)^{1/2}
\notag\\
&< \sigma_{k+1}(k+1)^{\alpha}\left(\int_k^{\infty}
\frac{1}{x^{2\alpha}} dx\right)^{1/2}  \mbox{\ \ due to $\alpha>\frac{1}{2}$}
\notag \\
&= \sigma_{k+1}\left(\frac{k+1}{k}\right)^{\alpha}
\sqrt{\frac{k}{2\alpha-1}}.
\label{modeest}
\end{align}
Substituting \eqref{modeest} into \eqref{error3} proves \eqref{error4}.
\qquad\endproof

From Theorems \ref{th1}--\ref{th2}, it is easy to see that
the error bounds \eqref{error5} and \eqref{error4} decrease with the
oversampling number $q$. Importantly,
whenever we take $q=k$ roughly, the factors in front of
$\sigma_{k+1}$ in \eqref{error5} and \eqref{error4}
reduce to $\mathcal{O}(1)$, independent of $n$, provided that
$\alpha$ is not close to $\frac{1}{2}$. On the other side,
for a fixed small $q\geq 4$, the factors in front of $\sigma_{k+1}$ in
\eqref{error5} and \eqref{error4} are $1+\mathcal{O}(\sqrt{k})$ and
$1+\mathcal{O}(k)$ for $\alpha$ not close to one, respectively, meaning
that the rank-$(k+q)$ approximation to $A$ may be more accurate for severely
ill-posed problems than for moderately and mildly ill-posed problems.
For a fixed $k$, the bigger $q$, the smaller
the bounds \eqref{error5} and \eqref{error4}, i.e., the
more accurate the rank-$k$ RSVD approximations.
As a result, in any event, our new bounds
are much sharper than \eqref{error2} and \eqref{overest},
and get more insight into the accuracy of rank-$(k+q)$ approximations
for $k$ not big and $\alpha$ not
close to $\frac{1}{2}$, where the factors in
front of $\sigma_{k+1}$ has been shown to lie between
$1+\mathcal{O}(\sqrt{n})$ and $1+\mathcal{O}(n)$.

Finally, we mention that all the results on $\|A-QQ^TA\|$ in this section
apply to $\|A-AQQ^T\|$ as well, where $AQQ^T$ is
generated by Algorithm~\ref{alg2}.

\section{TRSVD and error bounds, and the MTRSVD algorithms and their analysis}

We consider the MTRSVD method and compute the MTRSVD
solutions $x_{L,k}$ to the problem \eqref{mtik2} starting with
$k=1$. The MTRSVD solutions $x_{L,k}$ are regularized solutions to the general-form
regularization problem \eqref{general}.
MTRSVD first extracts a rank-$k$ TRSVD approximation
$\widetilde{A}_k$ from $QQ^TA$ (or $AQQ^T$) to $A$, and then utilizes
the LSQR algorithm \cite{paige1982} to iteratively solve
the resulting least squares problem at each iteration $k$ in \eqref{mtik2}.
This step is called inner iteration.
Starting with $k=1$, MTRSVD proceeds until a best regularized solution
is found at some $k=k_0$, at which the semi-convergence of MTRSVD occurs,
namely, the error $||L(x_{L,k}-x_{true})||$ decreases as $k$ increases up
to $k_0$ and then increases after $k>k_0$.

Recall that Algorithm \ref{alg1} generates a rank-$l$ approximation
$\widetilde{U}\widetilde{\Sigma}\widetilde{V}^T$ to $A$. Let
\begin{align*}
\widetilde{U} = (\widetilde{u}_1, \widetilde{u}_2, \ldots,\widetilde{u}_l)
\in \mathbb{R}^{m\times l},\quad \widetilde{V} = (\widetilde{v}_1, \widetilde{v}_2,
\ldots,\widetilde{v}_l)\in \mathbb{R}^{n\times l}
\end{align*}
and
$$\widetilde{\Sigma} = {\rm diag} (\widetilde{\sigma}_1, \widetilde{\sigma}_2, \ldots,
\widetilde{\sigma}_l)\in \mathbb{R}^{ l \times l}$$
with
$\widetilde{\sigma}_1\geq \widetilde{\sigma}_2\geq \cdots
\geq \widetilde{\sigma}_l>0$.
Take
\begin{align}\label{uk}
\widetilde{U}_k = (\widetilde{u}_1, \widetilde{u}_2, \ldots,\widetilde{u}_k)\in
\mathbb{R}^{m\times k},\quad
\widetilde{V}_k = (\widetilde{v}_1, \widetilde{v}_2, \ldots,\widetilde{v}_k)\in
\mathbb{R}^{n\times k},
\end{align}
and
\begin{align*}
\widetilde{\Sigma}_k = {\rm diag} (\widetilde{\sigma}_1, \widetilde{\sigma}_2, \ldots,
\widetilde{\sigma}_k)\in \mathbb{R}^{ k \times k},
\end{align*}
and form
\begin{align}\label{rankk}
\widetilde{A}_k=\widetilde{U}_k\widetilde{\Sigma}_k\widetilde{V}_k^T,
\end{align}
which is the best rank-$k$ approximation to
$\widetilde{U}\widetilde{\Sigma}\widetilde{V}^T=QQ^TA$,
called a rank-$k$ TRSVD approximation to $A$. Halko {\em et al}. \cite{Rev2011}
prove the following basic result.

\begin{theorem}
Let $\widetilde{A}_k$ be the rank-$k$ TRSVD approximation to $A$ defined by \eqref{rankk}.
Then the approximation error is
\begin{align}\label{error}
\|A-\widetilde{A}_k\|\leq \sigma_{k+1}+\|A-QQ^TA\|.
\end{align}
\end{theorem}

The bound \eqref{error} reflects the worse case. In Remark 9.1,
Halko {\em et al}. \cite{Rev2011} point out that
{\em ``In the randomized setting, the truncation step appears to be less
damaging than the error bound of Theorem 9.3 (i.e., \eqref{error} here)
suggests, but we currently lack a
complete theoretical understanding of its behavior."}
That is to say, the first
term $\sigma_{k+1}$ in \eqref{error} is generally conservative and
may be reduced substantially.

Keep in mind
that $A$ has $s$ singular values $\sigma_i$ with $s=\min\{m,n\}$.
Jia \cite{jia17c} has improved \eqref{error} and derived
sharper bounds, which explain why \eqref{error} may be an
overestimate, as shown in the following theorem.

\begin{theorem}[\cite{jia17c}]
Let $\widetilde{A}_k$ be the rank-$k$ TRSVD approximation to $A$ defined
by \eqref{rankk}. Then it holds that
\begin{align}\label{improv}
\|A-\widetilde{A}_k\|\leq \tilde{\sigma}_{k+1}+\|A-QQ^TA\|,
\end{align}
where $\tilde{\sigma}_{k+1}$ is the $(k+1)$-th singular value
of $Q^TA$ and satisfies
\begin{equation}\label{interrand}
\sigma_{m-q+1}\leq \tilde{\sigma}_{k+1}\leq \sigma_{k+1}
\end{equation}
with the definition $\sigma_{n+1}=\cdots=\sigma_m=0$. Analogously,
for the rank-$k$ TRSVD approximation $\widetilde{A}_k$
constructed by Algorithm~\ref{alg2}, it holds that
\begin{align}\label{improv2}
\|A-\widetilde{A}_k\|\leq \tilde{\sigma}_{k+1}+\|A-AQQ^T\|,
\end{align}
where $\tilde{\sigma}_{k+1}$ is the $(k+1)$-th singular value
of $AQ$ and satisfies
\begin{equation}\label{interrand2}
\sigma_{n-q+1}\leq \tilde{\sigma}_{k+1}\leq \sigma_{k+1}
\end{equation}
with the definition $\sigma_{m+1}=\cdots=\sigma_n=0$.

Particularly, the inequalities ``$\leq$" become strict ``$<$" in \eqref{interrand}
and \eqref{interrand2} if all the singular values $\sigma_j$ of $A$ are simple.
\end{theorem}

This theorem shows that the bound \eqref{improv} is unconditionally
superior to the bound
\eqref{error} and the former can improve the latter substantially
since $\tilde{\sigma}_{k+1}$ can be much smaller than $\sigma_{k+1}$ and
even be arbitrarily close to zero whenever $m-q+1>n$.
Once $\tilde{\sigma}_{k+1}<\sigma_{k+1}$ considerably,
the first term of \eqref{improv} is negligible relative to
the second term, and we will approximately have
\begin{align*}%
\|A-\widetilde{A}_k\|\approx \|A-QQ^TA\|.
\end{align*}

Regarding the MTRSVD solution $x_{L,k}$, we can establish the following result.

\begin{theorem}
Let $\widetilde{A}_k $ denote the rank-$k$ TRSVD approximation to $A$
generated by Algorithm~\ref{alg1} or Algorithm~\ref{alg2}.
Then the solution to \eqref{mtik2} can be written as
\begin{equation}\label{Xlksolu}
x_{L,k} = x_k-\left(L(I_n -\widetilde{V}_k\widetilde{V}_k^T)\right)^{\dag}Lx_k,
\end{equation}
where $x_k=\widetilde{A}_k^{+}b$ is the minimum 2-norm solution to the least
squares problem
\begin{align}\label{ls}
\min_{x\in\mathbb{R}^n}\|\widetilde{A}_kx-b\|.
\end{align}
\end{theorem}

{\em Proof}.
Following Eld${\rm \acute{e}}$n \cite{elden82}, we have
\begin{align}
x_{L,k}
&=(I_n-(L(I_n-\widetilde{A}_k^{\dagger}\widetilde{A}_k))^{+}L)\widetilde{A}_k^+b \nonumber\\
&=x_k-(L(I_n-\widetilde{A}_k^{\dagger}\widetilde{A}_k))^{+}Lx_k. \label{eq2}
\end{align}
Noting \eqref{rankk}, we have
$$\widetilde{A}_k^{\dagger}\widetilde{A}_k=\widetilde{V}_k\widetilde{V}_k^T.$$
Substituting the above into \eqref{eq2}, we obtain \eqref{Xlksolu}.
\qquad\endproof

Let $z_k = \left(L(I_n -\widetilde{V}_k\widetilde{V}_k^T)\right)^{\dag}Lx_k$.
Then $z_k$ is the minimum 2-norm solution to the least squares problem
\begin{equation}\label{eq3}
\min_{z\in \mathbb{R}^{n}}\|L(I_n -\widetilde{V}_k\widetilde{V}_k^T)z - Lx_k\|.
\end{equation}
We must point out that
the problem \eqref{eq3} and its coefficient matrix are different from those
in the MTSVD method in which the coefficient
matrix is $LV_{n-k}$ with $V_{n-k}=(v_{k+1},\ldots,v_n)$
available from the SVD \eqref{svd} of $A$.

Because of the large size of $L(I_n -\widetilde{V}_k\widetilde{V}_k^T)$,
we suppose that the problem \eqref{eq3} can only be solved by
iterative solvers. We will use
the LSQR algorithm \cite{paige1982} to solve the problem. In order
to make full use of the sparsity of $L$ itself and
reduce the computational cost and storage memory, it is vital to avoid
forming the dense matrix $L(I_n -\widetilde{V}_k\widetilde{V}_k^T)$ explicitly
within LSQR.
Notice that the only action of $L(I_n -\widetilde{V}_k\widetilde{V}_k^T)$
in the Lanczos diagonalization process and LSQR is to form
the products of it and its transpose with vectors.
We propose Algorithm~\ref{alg3}, which efficiently implements
the Lanczos bidiagonalization process with the starting vector
$\widehat{u}_1= Lx_k/\|Lx_k\|$.

\begin{algorithm}[htb]\label{alg3}
\caption{$\widehat k$-step Lanczos bidiagonalization process on
$L(I_n -\widetilde{V}_k\widetilde{V}_k^T)$}
\begin{algorithmic}[1]
\item[1.] Taking $\widehat{u}_1= Lx_k/\|Lx_k\|$, $w_1= L^T\widehat{u}_1$,
$g_1= \widetilde{V}_k^T\widehat{u}_1$ and define $\beta_1\widehat{v}_0 = 0$.
\item[2.] For $j=1, 2, \ldots, \widehat k$\\
$\widehat{p}=w_j-\widetilde{V}_k(\widetilde{V}_k^Tw_j)-\beta_j\widehat{v}_{j-1}$\\
$\alpha_j=\|\widehat{p}\|$;~~$\widehat{v}_j=\widehat{p}/\alpha_j$\\
$\widehat{r}=L\widehat{u}_j-L(\widetilde{V}_kg_j)-\alpha_j\widehat{u}_j$ \\
$\beta_{j+1} = \|\widehat{r}\|$;~~$\widehat{u}_{j+1} = \widehat{r}/\beta_{j+1}$\\
$w_{j+1} = L^T\widehat{u}_{j+1}$;~~$g_{j+1} = \widetilde{V}_k^T\widehat{u}_{j+1}$\\
\end{algorithmic}
\end{algorithm}

We now consider the solution of \eqref{eq3} using LSQR. Suppose
\begin{align}\label{vk}
\widetilde{V} =\left(
             \begin{array}{cc}
               \widetilde{V}_k & \widetilde{V}_{n-k} \\
             \end{array}
           \right)\in \mathbb{R}^{n\times n}
\end{align}
is an orthogonal matrix. It is then direct to obtain
\begin{align*}
L(I_n -\widetilde{V}_k\widetilde{V}_k^T)=L\widetilde{V}_{n-k}\widetilde{V}_{n-k}^T.
\end{align*}
Since $\widetilde{V}_{n-k}$ is column orthonormal, the nonzero singular values
of $L\widetilde{V}_{n-k}\widetilde{V}_{n-k}^T$ are identical to the singular values
of $L\widetilde{V}_{n-k}$. As a result, we have
\begin{align}\label{eq5}
\kappa(L(I_n -\widetilde{V}_k\widetilde{V}_k^T))&=\kappa(L\widetilde{V}_{n-k}
\widetilde{V}_{n-k}^T)=\kappa(L\widetilde{V}_{n-k}).
\end{align}

Next, we cite a lemma \cite[p.~78]{golub13} and exploit it to
investigate how the conditioning
of \eqref{eq3} changes as $k$ increases.

\begin{lemma}\label{lemma}
If
$B\in \mathbb{R}^{m\times n}$, $m>n$ and $c\in\mathbb{R}^m$, then
\begin{align*}
&\sigma_{\max}\left((
               \begin{array}{cc}
                 B & c \\
               \end{array}
             )\right)\geq \sigma_{\max}(B),\\
&\sigma_{\min}\left((
               \begin{array}{cc}
                 B & c \\
               \end{array}
             )\right)\leq \sigma_{\min}(B).
\end{align*}
\end{lemma}

This lemma shows that if a column is added to a rectangular matrix
then the largest singular value increases and the smallest singular value
decreases. Therefore, we directly obtain the following result on the
conditioning of \eqref{eq3}.

\begin{theorem}\label{thm5}
Let the matrix $\widetilde{V}_{n-k}$ be defined by \eqref{vk}. Then for
$p\geq n-k$, we have
\begin{align}\label{condition}
\kappa(L\widetilde{V}_{n-k})\geq\kappa(L\widetilde{V}_{n-(k+1)}),\quad k = 1, 2,
\ldots, n-1,
\end{align}
i.e.,
\begin{align}\label{condition2}
\kappa(L(I_n -\widetilde{V}_k\widetilde{V}_k^T))\geq
\kappa(L(I_n -\widetilde{V}_{k-1}\widetilde{V}_{k-1}^T)),\quad k =1, 2,
\ldots,n-1.
\end{align}
\end{theorem}

This theorem indicates that, when applied to solving \eqref{eq3},
the LSQR algorithm generally converges faster with $k$
by recalling that the worst convergence factor of LSQR
is $\frac{\kappa(L\widetilde{V}_{n-k})+1}{\kappa(L\widetilde{V}_{n-k})-1}$;
see \cite[p.~291]{bjorck96}.
Particularly, in exact arithmetic, LSQR will find the exact solution
$z_k$ of \eqref{eq3} after at most $n-k$ iterations.

Having done the above, we can present our MTRSVD algorithm for
the $m\geq n$ case, named as Algorithm \ref{alg4}.

\begin{algorithm}[htb]
\caption{(MTRSVD) Given $A\in \mathbb{R}^{m\times n}(m\geq n)$ and $l=k+q<n$
and $q\geq4$, compute the solution $x_{L,k}$ of \eqref{mtik2}.}

\begin{algorithmic}[1]\label{alg4}

\STATE Use Algorithm~\ref{alg1} to compute the rank-$k$ TRSVD
approximation $\widetilde{A}_k$ to $A$:
$\widetilde{A}_k=\widetilde{U}_k\widetilde{\Sigma}_k\widetilde{V}_k^T$.

\STATE Compute the the minimum 2-norm solution $x_k$ to \eqref{ls}.

\STATE Compute the solution $z_k$ to \eqref{eq3} by LSQR.

\STATE Compute the solution $x_{L,k}$, defined by \eqref{Xlksolu},
to the problem \eqref{mtik2}.
\end{algorithmic}
\end{algorithm}

For the $m\leq n$ case, making use of Algorithm \ref{alg2},
we present Algorithm \ref{alg5}, a variant of Algorithm \ref{alg4}.

\begin{algorithm}[htb]
\caption{(MTRSVD) Given $A\in \mathbb{R}^{m\times n}(m\leq n)$ and $l=k+q<n$ and
$q\geq4$, compute the solution $x_{L,k}$ of \eqref{mtik2}.}
\begin{algorithmic}[1]\label{alg5}
\STATE Use Algorithm~\ref{alg2} to compute the rank-$k$ TRSVD approximation
$\widetilde{A}_k$ to $A$:
$\widetilde{A}_k= \widetilde{U}_k\widetilde{\Sigma}_k\widetilde{V}_k^T$.

\STATE Compute the minimum 2-norm solution $x_k$ to \eqref{ls}.

\STATE Compute the solution $z_k$ to \eqref{eq3} by LSQR.

\STATE Compute the solution $x_{L,k}$, defined by \eqref{Xlksolu}, to the problem
\eqref{mtik2}.
\end{algorithmic}
\end{algorithm}

We comment that at step 3 of Algorithms~\ref{alg4}--\ref{alg5}, in numerical
experiments we will use the Matlab function {\sf lsqr.m} to solve the problems with
a given tolerance $tol$ as the stopping
criterion. In what follows we make a detailed analysis and
show that the default $tol=10^{-6}$ is generally
good enough  and larger $tol$ can be allowed in practical applications.

First of all, let us estimate the accuracy of the computed solution $\bar{z}_k$
with the stopping
tolerance $tol$.
Let $r= Lx_k-L(I_n -\widetilde{V}_k\widetilde{V}_k^T)z_k$ be
the residual of the solution $z_k$ to the problem \eqref{eq3}.
It is known from \cite{paige1982} that, with the
stopping tolerance $tol$, the computed $\bar{z}_k$ is the exact
solution to the perturbed problem
\begin{equation}\label{perb}
\min_{z\in \mathbb{R}^{n}}\|(L(I_n -\widetilde{V}_k\widetilde{V}_k^T)+E_k)z - Lx_k\|,
\end{equation}
where the perturbation matrix
$$
E_k=-\frac{r_kr_k^TL(I_n -\widetilde{V}_k\widetilde{V}_k^T)}{\|r_k\|^2}.
$$
with
$$
r_k=Lx_k
-L(I_n -\widetilde{V}_k\widetilde{V}_k^T)\bar{z}_k
$$
being the residual of the computed solution $\bar{z}_k$ and
$$
\frac{\|E_k\|}{\|L(I_n -\widetilde{V}_k\widetilde{V}_k^T)\|}
=\frac{\|(I_n -\widetilde{V}_k\widetilde{V}_k^T)L^Tr_k\|}
{\|L(I_n -\widetilde{V}_k\widetilde{V}_k^T)\|\|r_k\|}\leq tol.
$$
For details on implementations, we refer to \cite{paige1982}.

With the above notation and \eqref{eq5}, defining $\eta=tol\cdot
\kappa(L\widetilde{V}_{n-k})$, exploiting
the standard perturbation theory \cite[p.~382]{higham02}, we obtain
\begin{equation}\label{relerror}
\frac{\|z_k-\bar{z}_k\|}{\|z_k\|}\leq
\frac{tol\cdot \kappa(L\widetilde{V}_{n-k})}
{1-\eta}\left(2+(\kappa(L\widetilde{V}_{n-k})+1)\frac{\|r\|}
{\|L(I_n -\widetilde{V}_k\widetilde{V}_k^T)\|\|z_k\|}
\right).
\end{equation}
Actually, by checking its proof we find that the above factor
$\kappa(L\widetilde{V}_{n-k})+1$
can be replaced by $\kappa(L\widetilde{V}_{n-k})$
in our context since the left hand side $Lx_k$ in the perturbed \eqref{perb}
is unperturbed.

In applications, $L$ is typically well conditioned
\cite{hansen98,hansen10}.
Since $\kappa(L)\geq \kappa(L\widetilde{V}_{n-k})$ for $p\geq n-k$,
the left hand side of \eqref{relerror} is {\em at least} as small as ${\cal O}(tol)$
with a generic constant in ${\cal O}(\cdot)$.

Recall from \eqref{Xlksolu} that the MTRSVD solution
$$
x_{L,k}=x_k-z_k,
$$
and define the computed solution
$$
\bar{x}_{L,k}=x_k-\bar{z}_k.
$$
We thus have $\|x_{L,k}-\bar{x}_{L,k}\|=\|z_k-\bar{z}_k\|$,
from which and \eqref{relerror}
it is reasonable to suppose
\begin{equation}\label{xkzk}
\frac{\|x_{L,k}-\bar{x}_{L,k}\|}{\|x_{L,k}\|}
\approx\frac{\|z_k-\bar{z}_k\|}{\|z_k\|}\leq {\cal O}(tol)
\end{equation}
since it is generally impossible that
$\|x_{L,k}\|$ is much smaller or larger than $\|z_k\|$.

Let $x_{L}^{opt}$ be a {\em best possible} regularized solution to
the problem \eqref{general} with the white noise $e$.
Then under a certain necessary discrete
Picard condition, a GSVD analysis indicates
that the error $\|x_{L}^{opt}-x_{true}\|\geq {\cal O}(\|e\|)$ with a
generic constant in ${\cal O}(\cdot)$; see \cite[p.~83]{hansen98}. 

Let $x_{L,k_0}$ be the best regularized solutions by the MTRSVD algorithms. Then
\begin{equation}\label{compare}
\|x_{L,k_0}-x_{true}\|\geq \|x_L^{opt}-x_{true}\|\geq {\cal O}(\|e\|).
\end{equation}
By \eqref{relerror} and \eqref{xkzk} as well as $\|x_{L,k_0}\|\approx \|x_{true}\|$,
we have
\begin{eqnarray}
\frac{\|\bar{x}_{L,k_0}-x_{true}\|}{\|x_{true}\|}&\leq &\frac{\|x_{L,k_0}-x_{true}\|}
{\|x_{true}\|}+\frac{\|x_{L,k_0}-\bar{x}_{L,k_0}\|}{\|x_{true}\|} \nonumber\\
&=&\frac{\|x_{L,k_0}-x_{true}\|}
{\|x_{true}\|}+\frac{\|x_{L,k_0}-\bar{x}_{L,k_0}\|}{\|x_{L,k_0}\|}
\frac{\|x_{L,k_0}\|}{\|x_{true}\|} \nonumber\\
&=&\frac{\|x_{L,k_0}-x_{true}\|}
{\|x_{true}\|}+{\cal O}(tol). \nonumber
\end{eqnarray}
On the other hand, we similarly obtain
$$
\frac{\|\bar{x}_{L,k_0}-x_{true}\|}{\|x_{true}\|}\geq \frac{\|x_{L,k_0}-x_{true}\|}
{\|x_{true}\|}-{\cal O}(tol).
$$

Suppose that the noise free problem of \eqref{eq1} is consistent,
i.e., $Ax_{true}=b_{true}$.
Since $\|A\|\|x_{true}\|\geq \|b_{true}\|$, it follows from \eqref{compare} that
$$
\frac{\|x_{L,k_0}-x_{true}\|}
{\|x_{true}\|}\approx \|A\|\frac{\|x_{L,k_0}-x_{true}\|}
{\|b_{true}\|}\geq \|A\|{\cal O}\left(\frac{\|e\|}{\|b_{true}\|}\right)
={\cal O}\left(\frac{\|e\|}{\|b_{true}\|}\right)
$$
when $\|A\|\approx 1$
(this can always be done by suitable
scaling). As a result, summarizing the above derivation, we have proved the following
results.

\begin{theorem}\label{thm6}
If $L$ is well conditioned and
\begin{equation}\label{accucond}
{\cal O}(tol)< \frac{\|e\|}{\|b_{true}\|},
\end{equation}
then
\begin{equation}\label{appeq}
\frac{\|x_{L,k_0}-x_{true}\|}
{\|x_{true}\|}-{\cal O}(tol)
\leq \frac{\|\bar{x}_{L,k_0}-x_{true}\|}{\|x_{true}\|}\leq
\frac{\|x_{L,k_0}-x_{true}\|}{\|x_{true}\|}+{\cal O}(tol),
\end{equation}
i.e.,
\begin{equation}\label{computerror}
\frac{\|\bar{x}_{L,k_0}-x_{true}\|}{\|x_{true}\|}=\frac{\|x_{L,k_0}-x_{true}\|}
{\|x_{true}\|}
\end{equation}
within the error ${\cal O}(tol)$ with a generic constant in ${\cal O}(\cdot)$,
meaning that the computed $\bar{x}_{L,k_0}$ has the same
as the exact $x_{L,k_0}$ as an approximation to $x_{true}$.
\end{theorem}

Furthermore, based the above, we can establish general results, which
include Theorem~\ref{thm6} as a special case.
Since $x_{L,k_0}$'s are best possible regularized solutions
by the MTRSVD algorithms, i.e.,
$$
\frac{\|x_{L,k_0}-x_{true}\|}
{\|x_{true}\|}=\min_{k=1,2,\ldots,n} \frac{\|x_{L,k}-x_{true}\|}
{\|x_{true}\|},
$$
under the condition \eqref{accucond}, it
follows from the fact
\begin{equation}\label{estimate}
\frac{\|x_{L,k}\|}{\|x_{true}\|}={\cal O}(1)
\end{equation}
and the proof of Theorem~\ref{thm6} that
\eqref{appeq} and \eqref{computerror} also hold
when the index $k_0$ is replaced by $k=1,2,\ldots,k_0$ and {\em a few} $k>k_0$.
We remark that the estimate \eqref{estimate} holds because $\|x_{L,k}\|$
exhibits increasing tendency, and it first
approximates $\|x_{true}\|$ from below for $k=1,2,\ldots,k_0$ and
then starts to deviate from $\|x_{true}\|$ but not too much for a few
$k>k_0$. Therefore, we have proved the following theorem.

\begin{theorem}\label{thm7}
If $L$ is well conditioned and
$$
{\cal O}(tol)< \frac{\|e\|}{\|b_{true}\|},
$$
then for $k=1,2,\ldots k_0$ and a few $k>k_0$ we have
\begin{equation}\label{appeq2}
\frac{\|x_{L,k}-x_{true}\|}
{\|x_{true}\|}-{\cal O}(tol)
\leq \frac{\|\bar{x}_{L,k}-x_{true}\|}{\|x_{true}\|}\leq
\frac{\|x_{L,k}-x_{true}\|}{\|x_{true}\|}+{\cal O}(tol),
\end{equation}
i.e.,
\begin{equation}\label{computerror2}
\frac{\|\bar{x}_{L,k}-x_{true}\|}{\|x_{true}\|}=\frac{\|x_{L,k}-x_{true}\|}
{\|x_{true}\|}
\end{equation}
within the error ${\cal O}(tol)$ with a generic constant in ${\cal O}(\cdot)$,
meaning that the computed $\bar{x}_{L,k}$ has the same
as the exact $x_{L,k}$ as an approximation to $x_{true}$.
\end{theorem}

It is worthwhile to notice that
the relative noise level $\frac{\|e\|}{\|b_{true}\|}$ is typically
more or less around $10^{-3}$ in applications, three orders bigger than
$10^{-6}$.
Combining all the above together,
we come to conclude that it is generally enough to set $tol=10^{-6}$ in LSQR at step
3 of Algorithms~\ref{alg4}--\ref{alg5}. A smaller $tol$ will result in more inner
iterations without any gain in the accuracy of $\bar{x}_{L,k}$ as
regularized solutions for $k=1,2,\ldots,k_0$ and a few $k>k_0$.
Moreover, Theorems~\ref{thm6}--\ref{thm7} indicate that $tol=10^{-6}$
is generally well conservative and larger $tol$ can be used, so
that LSQR uses fewer iterations to achieve the convergence and the MTRSVD algorithms
are more efficient.

In summary, our conclusion is that a widely varying choice of $tol$ has {\em no}
effects of regularization of the MTRSVD algorithms, provided
that $tol<\frac{\|e\|}{\|b_{true}\|}$ considerably and
the regularization matrix $L$ is well conditioned, but it
has {\em substantial effects} on the efficiency of MTRSVD.
In our numerical experiments, we have found that for each test problem
with $\frac{\|e\|}{\|b_{true}\|}=10^{-2}$ and $10^{-3}$
the computed best regularized solutions obtained
by MTRSVD have the same accuracy and the
convergence curves of MTRSVD are indistinguishable when
taking three $tol=10^{-6}, 10^{-5}$ and $10^{-4}$.

For a given oversampling parameter $q\geq 4$, we need to determine an optimal $k=k_0$
for finding a best possible regularized solution $x_{L,k_0}$ in the MTRSVD
algorithms. It is crucial to realize that,
just as in the TSVD method, the parameter $k$ plays
the role of the regularization parameter in the MTSVD, TGSVD, TRGSVD
and MTRSVD methods. From now on,
denote by $x_k^{reg}$ the regularized solution at step $k$
obtained by each of them. These methods must exhibit semi-convergence
\cite{hansen98,hansen10,nat1986}:
the error $\|L(x_k^{reg}-x_{true})\|$ decreases (correspondingly,
$\|Lx_k^{reg}\|$ steadily increases) with respect to $k$ in
the first stage until some step $k=k_0$ and then starts to increases
(correspondingly, $\|Lx_k^{reg}\|$
starts to increase considerably) after $k>k_0$. Such $k_0$
is exactly an optimal regularization parameter, at which
the regularized solution $x_{k_0}^{reg}$
is most accurate and is thus the best possible one obtained by each of these
methods.

Given an oversampling parameter $q$, the algorithms of Wei {\em et al}.
\cite{wei16} and Xiang and Zou \cite{zou15} first generate RGSVDs for a
certain fixed $k$ and then determine the optimal Tikhonov regularization
parameter $\lambda_{opt}$ by GCV \cite{hansen98,hansen10}. In the $m\geq n$ case,
Wei {\em et al}. \cite{wei16} use Algorithm 4.2
in \cite{Rev2011} to determine such a $k$ adaptively until
\begin{align*}
\|A-QQ^TA\|\leq\widetilde{\varepsilon}
\end{align*}
is satisfied for some small $\widetilde{\varepsilon}$. Then they
replace $A$ by the truncated rank-$k$ approximation to $A$ obtained from
the GSVD of $\{Q^TA,L\}$ in
\eqref{general3} and determine $\lambda_{opt}$; see (2.8) in \cite{wei16}.
For the $m\leq n$ case, they reduce the original large \eqref{general3} to
a projected problem that replaces $A$ and $L$ by $AQ$ and $LQ$, respectively,
with $\|A-AQQ^T\|\leq \widetilde{\varepsilon}$, and then solve it by using
the GSVD of $\{AQ,LQ\}$ and determining an optimal $\lambda_{opt}$.
In the numerical experiments, they take
a fixed $\widetilde{\varepsilon}=10^{-2}$ for all the test problems and
noise levels. They emphasize that the choice of an optimal tolerance
$\widetilde{\varepsilon}$ is an open problem. As a matter of fact,
the size of $\widetilde{\varepsilon}$ is
problem and noise level dependent, and it is impossible to
presume a fixed and optimal $\widetilde{\varepsilon}$
for all problems and noise levels.
A basic fact is that the smaller the noise level, the more dominant SVD (or GSVD)
components of $A$ are needed \cite{hansen98,hansen10} to form
best regularized solutions. This means that
the smaller the noise level, the smaller $\widetilde{\varepsilon}$ must be.

In practical applications, for the TSVD and MTSVD methods, one can use
the GCV parameter-choice method or the L-curve criterion to determine their
regularization parameters $k_0$ \cite{hansen98,hansen10,hansen92}.
The L-curve criterion is directly applicable to our MTRSVD algorithms:
Given an oversampling parameter $q$, they proceed from $k=1$ onwards,
successively increment $l=k+q$ and expand $Q$. The algorithms compute a sequence
of regularized solutions $x_k^{reg}$, and we plot the curve of
$(||Ax_k^{reg}-b\|,\|Lx_k^{reg}\|)$
in log-log scale, whose corner corresponds to the best regularized
solution $x_{k_0}^{reg}$ with $k_0$ the optimal regularization parameter,
at which the semi-convergence of our algorithms occurs.
In contrast, the GCV parameter-choice method
is not directly applicable to our MTRSVD algorithms,
and some nontrivial effects are needed to derive corresponding GCV
functions. We will consider
the GCV parameter-choice method for our MTRSVD algorithms in future work.

In our next numerical
experiments, the true solutions $x_{true}$'s to all the test
problems are known, so that for a sequence of regularized solutions
$x_k^{reg}$ the a-priori relative errors $\|L(x_k^{reg}-x_{true})\|/
\|Lx_{true}\|$ can be computed and the regularization parameter
$k_0$ of semi-convergence is easily identified for each method
by plotting the corresponding convergence curve.

\section{Numerical examples}

In this section, we report numerical experiments to demonstrate that the MTRSVD
algorithms can compute regularized solutions as accurately as
the standard TGSVD algorithm and at least as accurately as those obtained by
the RGSVD algorithms in \cite{wei16} and \cite{zou15}.
We choose some one dimensional examples from Hansen's
regularization toolboxs \cite{hansen2007(2)} and a two dimensional problem
from \cite{hansen12}.
We generated the Gaussian noise vectors $e$ whose entries are normally distributed
with mean zero. We denote the relative noise level $\varepsilon = \frac{\|e\|}{\|b_{true}\|}$,
and use $\varepsilon = 10^{-2}, 10^{-3}$ in the experiments. To simulate exact arithmetic,
the full reorthogonalization is used during the Lanczos bidiagonalization process.
Purely for test purposes, we choose $L=L_1$ and $L_3$ defined by \eqref{l1}
and \eqref{l1l2}, respectively. For $L=L_2$, we have found that the results and
comparisons are very similar
to those for $L=L_1$, so we omit the reports on $L=L_2$.

Recall that $x_k^{reg}$ denotes the regularized solution obtained by
each of TGSVD, MTRSVD and RGSVD. We use the the relative error
$$
\frac{\|L(x_k^{reg}-x_{true})\|}{\|Lx_{true}\|}
$$
to plot the convergence curve of each method with respect to $k$.
The TRGSVD algorithms in \cite{wei16} and \cite{zou15} are
denoted by {\sf weirgsvd} and {\sf xiangrgsvd}, respectively,
we abbreviate the standard TGSVD algorithm as {\sf tgsvd} and
Algorithms~\ref{alg4}--\ref{alg5} as {\sf mtrsvd}. Here
we make some non-essential modifications on the original
{\sf weirgsvd} and {\sf xiangrgsvd} in order to
compare all the algorithms under consideration more directly and insightfully.
The original RGSVD algorithms in \cite{wei16, zou15} are
the combinations of RGSVD and general-form
Tikhonov regularization. We now truncate
rank-$(k+q)$ RGSVD and obtain a rank-$k$ truncated randomized GSVD (TRGSVD),
leading to the corresponding TRGSVD algorithms, such a TRGSVD algorithm
was mentioned by Xiang and Zou \cite{zou15}. The original
{\sf weirgsvd} and {\sf xiangrgsvd} and the current ones are the same
in the spirit of {\sf tgsvd} and the GSVD with
Tikhonov regularization \cite{hansen98,hansen10}, and they will
generate the best regularized solutions
with essentially the same accuracy. In the tables to
be presented, we will list the given oversampling parameter $q$ and the optimal
regularization parameter $k_0$ in the braces.  We use the
Matlab function {\sf lsqr.m} to solve the least squares problems \eqref{eq3}
with the default stopping tolerance $tol=10^{-6}$. We have observed that for
$\varepsilon=10^{-2},10^{-3}$ the three convergence curves of {\sf mtrsvd} are
indistinguishable for each test problem when taking $tol=10^{-6},10^{-5},10^{-4}$,
respectively, and the computed best regularized solutions by {\sf mtrsvd} for these
three $tol$ have the same accuracy. As a result, for the sake of uniqueness
and length, we will only report the results on $tol=10^{-6}$.

All the computations are carried out in Matlab R2015b 64-bit on
Intel Core i3-2120 CPU 3.30GHz processor and 4 GB RAM.

\subsection{The $m\geq n$ case}

We first present the results on four one dimensional test problems from Hansen's
regularization toolbox \cite{hansen2007(2)}, and then report the results on
a two dimensional test problem from Hansen's regularization toolbox \cite{hansen12}.

\subsubsection{The one dimensional case}

All test problems arises from the discretization of the
first kind Fredholm integral equations
\begin{equation}\label{integral}
\int_a^b k(s,t)x(t)dt = f(s), \quad c\leq s\leq d.
\end{equation}
For each problem we use the code of \cite{hansen2007(2)} to generate $A$,
the true solution $x_{true}$ and noise-free right-hand side $b_{true}$.
The four test problems are severely, moderately and mildly ill-posed, respectively;
see Table \ref{tab1}, where we choose
the parameter "$example=2$" for the test problem {\sf deriv2}.

\begin{table}[htp]
    \centering
    \caption{The description of test problems.}
    \begin{tabular}{lll}
     \hline
     Problem  & Description                                & Ill-posedness \\
     \hline
     {\sf shaw}     &    One dimensional image restoration model & severe\\
     {\sf gravity}  & One dimensional gravity surveying problem  & severe\\
     {\sf heat}     & Inverse heat equation                      & moderate\\
     {\sf deriv2}   & Computation of second derivative           & mild\\
     \hline
   \end{tabular}
   \label{tab1}
\end{table}

\begin{table}[htp]
  \centering
  \caption{The comparison of Algorithm \ref{alg4} ({\sf mtrsvd})
  and the others with $L=L_1$ and
  $\varepsilon=10^{-2}, \ 10^{-3}$.}
   \label{tab2}
\centerline{$\varepsilon=10^{-2}$}
     \begin{minipage}[t]{1\textwidth}
     \begin{tabular*}{\linewidth}{lp{0.5cm}p{1.7cm}p{1.7cm}p{1.7cm}p{1.7cm}p{1.7cm}p{1.7cm}p{1.7cm}}
     \toprule[0.6pt]
     \multicolumn{2}{c}{}   & \multicolumn{3}{c}{$m=n=1,024$}&\multicolumn{2}{c}{$m=n=10,240$} \\
     \cmidrule(lr){3-5}\cmidrule(lr){6-7}
                     &$q$      &{\sf tgsvd}  &{\sf weirgsvd}&{\sf mtrsvd}
                     	    &{\sf weirgsvd}      &{\sf mtrsvd}  \\ \midrule[0.6pt]
     {\sf shaw}      &9      & 0.2043(6)    &0.2043(7)	  &0.2043(7)	& \quad\quad-   &0.1946(7)\\
     {\sf gravity}   &11      & 0.3205(7)	&0.3203(8)	  &0.3202(8)	& \quad\quad-   &0.2594(9)\\
     {\sf heat}      &7      & 0.2526(23)	&0.2544(23)	  &0.2457(23)   & \quad\quad-   &0.2285(23)\\
     {\sf deriv2}    &11      & 0.4264(5)	&0.4324(6)	  &0.4411(6)	& \quad\quad-   &0.3621(16)\\
    \bottomrule[0.6pt]
     \end{tabular*}\\[2pt]
  \end{minipage}
\centerline{$\varepsilon=10^{-3}$}
  \begin{minipage}[t]{1\textwidth}
     \begin{tabular*}{\linewidth}{lp{0.5cm}p{1.7cm}p{1.7cm}p{1.7cm}p{1.7cm}p{1.7cm}p{1.7cm}p{1.7cm}}
     \toprule[0.6pt]
     \multicolumn{2}{c}{}   &\multicolumn{3}{c}{$m=n=1,024$} &\multicolumn{2}{c}{$m=n=10,240$} \\
     \cmidrule(lr){3-5}\cmidrule(lr){6-7}
                      &$q$   &{\sf tgsvd}        &{\sf weirgsvd}	    &{\sf mtrsvd}
                      &{\sf weirgsvd}  &{\sf mtrsvd}  \\ \midrule[0.6pt]
     {\sf shaw}      &9      & 0.1681(8)    &0.1681(9)	  &0.1681(9)	& \quad\quad-   &0.1428(9)\\
     {\sf gravity}   &7      & 0.2660(10)	&0.2675(11)	  &0.2660(11)	& \quad\quad-   &0.2532(11)\\
     {\sf heat}      &8      & 0.1664(30)	&0.1673(31)	  &0.1623(29)   & \quad\quad-   &0.1399(36)\\
     {\sf deriv2}    &6      & 0.3341(11)	&0.3360(12)	  &0.3462(12)	& \quad\quad-   &0.2916(12)\\
    \bottomrule[0.6pt]
     \end{tabular*}\\[2pt]
  \end{minipage}
\end{table}

In Table~\ref{tab2}, we display the relative errors of
the best regularized solutions $x_{k_0}^{reg}$
by {\sf tgsvd}, {\sf weirgsvd} and {\sf mtrsvd}
with $L = L_1$ and $\varepsilon =10^{-2}, 10^{-3}$,
respectively. They illustrate that for all test problems
with $m=n=1,024$ the solution accuracy of {\sf mtrsvd}
is very comparable to that of {\sf tgsvd} and {\sf weirgsvd}.
For $m=n=10,240$, {\sf tgsvd} and {\sf weirgsvd}
are out of memory in our computer,
but {\sf mtrsvd} works well and the best
regularized solution is more accurate than the corresponding
one for $m=n=1,024$. We observe from the table
that for each test problem the best regularized solution by each algorithm
is correspondingly more accurate for $\varepsilon =10^{-3}$ than
$\varepsilon =10^{-2}$; for each algorithm,
the optimal regularization parameter $k_0$ is bigger for $\varepsilon =10^{-3}$
than for $\varepsilon =10^{-2}$. All these are expected and justify
that the smaller the noise level $\varepsilon$ is, the more SVD (or GSVD) dominant
components of $A$ or ($\{A,L\}$) are needed to form best regularized solutions.
Finally, as is seen, for each problem and the given $\varepsilon$,
the optimal $k_0$ are almost the same for all the algorithms. This indicates
that {\sf mtrsvd} and {\sf weirgsvd} effectively capture the dominant SVD and GSVD
components of $A$ and $\{A,L\}$, respectively.

\begin{table}[htp]
  \centering
  \caption{The comparison of Algorithm \ref{alg4} ({\sf mtrsvd}) and the others with
  $L=L_3$ and $\varepsilon=10^{-2}, \ 10^{-3}$.}
    \label{tab4}
    \centerline{$\varepsilon=10^{-2}$}
     \begin{minipage}[t]{1\textwidth}
     \begin{tabular*}{\linewidth}{lp{0.5cm}p{1.7cm}p{1.7cm}p{1.7cm}p{1.7cm}p{1.7cm}p{1.7cm}p{1.7cm}}
     \toprule[0.6pt]
     \multicolumn{2}{c}{}   &\multicolumn{3}{c}{$m=n=1,024$} &\multicolumn{2}{c}{$m=n=10,240$} \\
     \cmidrule(lr){3-5}\cmidrule(lr){6-7}
                     &$q$   &{\sf tgsvd}  &{\sf weirgsvd} &{\sf mtrsvd}	 &{\sf weirgsvd} &{\sf mtrsvd}
                      \\ \midrule[0.6pt]
     {\sf shaw}      &11   &0.2030(7)	&0.2030(8)	   &0.2024(8)	& \quad\quad-        &0.1984(7)\\
     {\sf gravity}   &10   &0.3340(8)	&0.3342(8)	   &0.3339(8)	& \quad\quad-        &0.2292(9)\\
     {\sf heat}      &7   &0.2966(23)	&0.2856(23)	   &0.2695(22)	& \quad\quad-        &0.2386(23)\\
     {\sf deriv2}    &9   &0.4365(6)	&0.4446(6)     &0.4430(7)	& \quad\quad-        &0.4207(10)\\
     \bottomrule[0.6pt]
     \end{tabular*}\\[2pt]
  \end{minipage}
   \centerline{$\varepsilon=10^{-3}$}
   \begin{minipage}[t]{1\textwidth}
     \begin{tabular*}{\linewidth}{lp{0.5cm}p{1.7cm}p{1.7cm}p{1.7cm}p{1.7cm}p{1.7cm}p{1.7cm}p{1.7cm}}
     \toprule[0.6pt]
     \multicolumn{2}{c}{}   &\multicolumn{3}{c}{$m=n=1,024$} &\multicolumn{2}{c}{$m=n=10,240$} \\
     \cmidrule(lr){3-5}\cmidrule(lr){6-7}
                     &$q$   &{\sf tgsvd} &{\sf weirgsvd} &{\sf mtrsvd}&{\sf weirgsvd} &{\sf mtrsvd}
                     \\ \midrule[0.6pt]
     {\sf shaw}      &4   &0.1694(8)	&0.1694(8)	   &0.1694(8)	& \quad\quad-        &0.1431(9)\\
     {\sf gravity}   &8   &0.2838(10)	&0.2830(10)	   &0.2811(10)	& \quad\quad-        &0.1789(9)\\
     {\sf heat}      &12   &0.1626(30)	&0.1616(30)	   &0.1610(30)	& \quad\quad-        &0.1468(35)\\
     {\sf deriv2}    &8   &0.3465(10)	&0.3499(10)    &0.3550(10)	& \quad\quad-        &0.3129(13)\\
      \bottomrule[0.6pt]
     \end{tabular*}\\[2pt]
  \end{minipage}
\end{table}

In Table~\ref{tab4}, we display the relative
errors of the best regularized solutions by {\sf tgsvd}, {\sf weirgsvd} and {\sf mtrsvd}
with $L = L_3$ and $\varepsilon =10^{-2}, 10^{-3}$,
respectively. The results and performance evaluations
on the three algorithms are analogous to those for $L=L_1$, and the details
are thus omitted.

\begin{figure}
\begin{minipage}{0.48\linewidth}
  \centerline{\includegraphics[width=6.0cm,height=4cm]{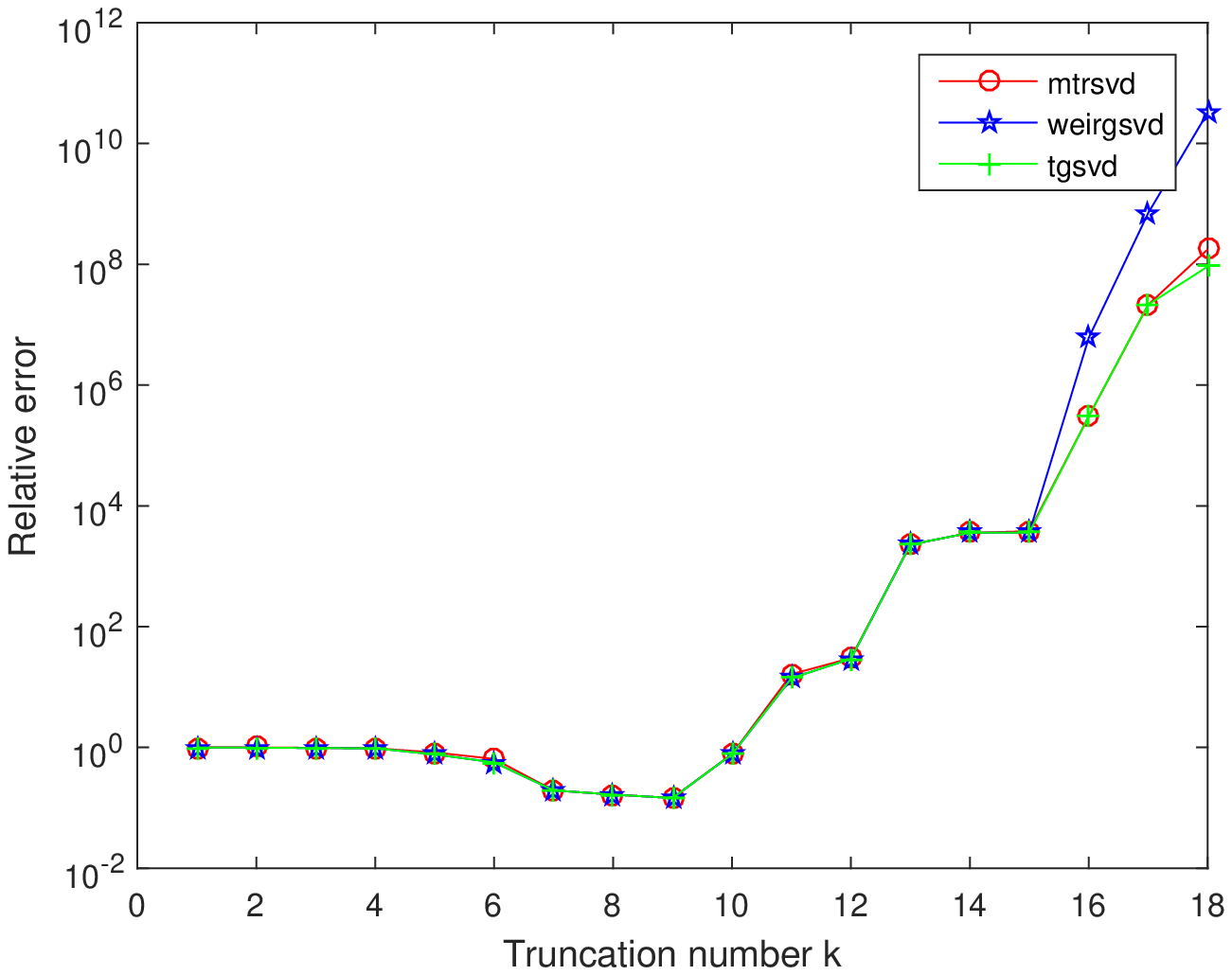}}
  \centerline{(a)}
\end{minipage}
\hfill
\begin{minipage}{0.48\linewidth}
  \centerline{\includegraphics[width=6.0cm,height=4cm]{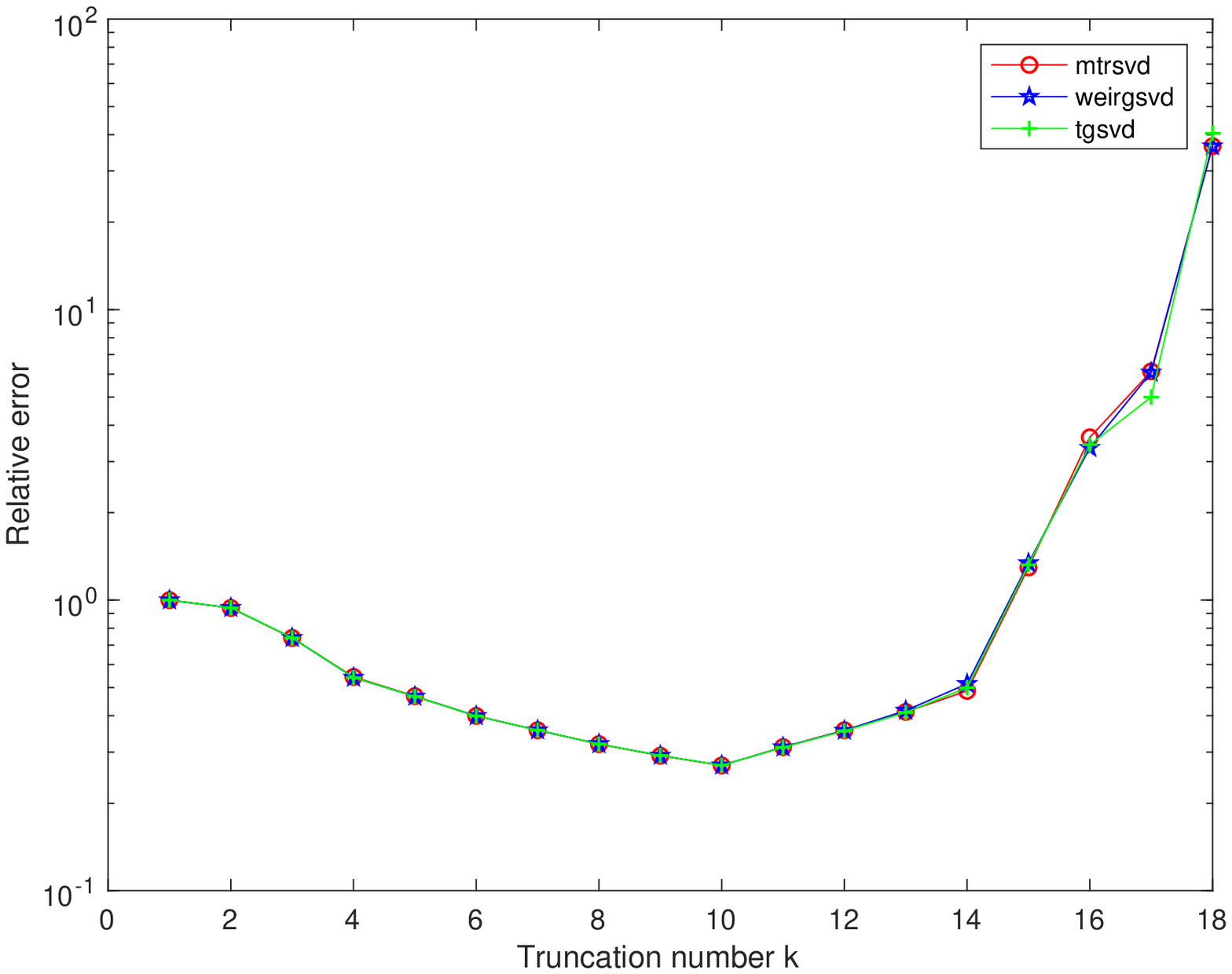}}
  \centerline{(b)}
\end{minipage}

\vfill
\begin{minipage}{0.48\linewidth}
  \centerline{\includegraphics[width=6.0cm,height=4cm]{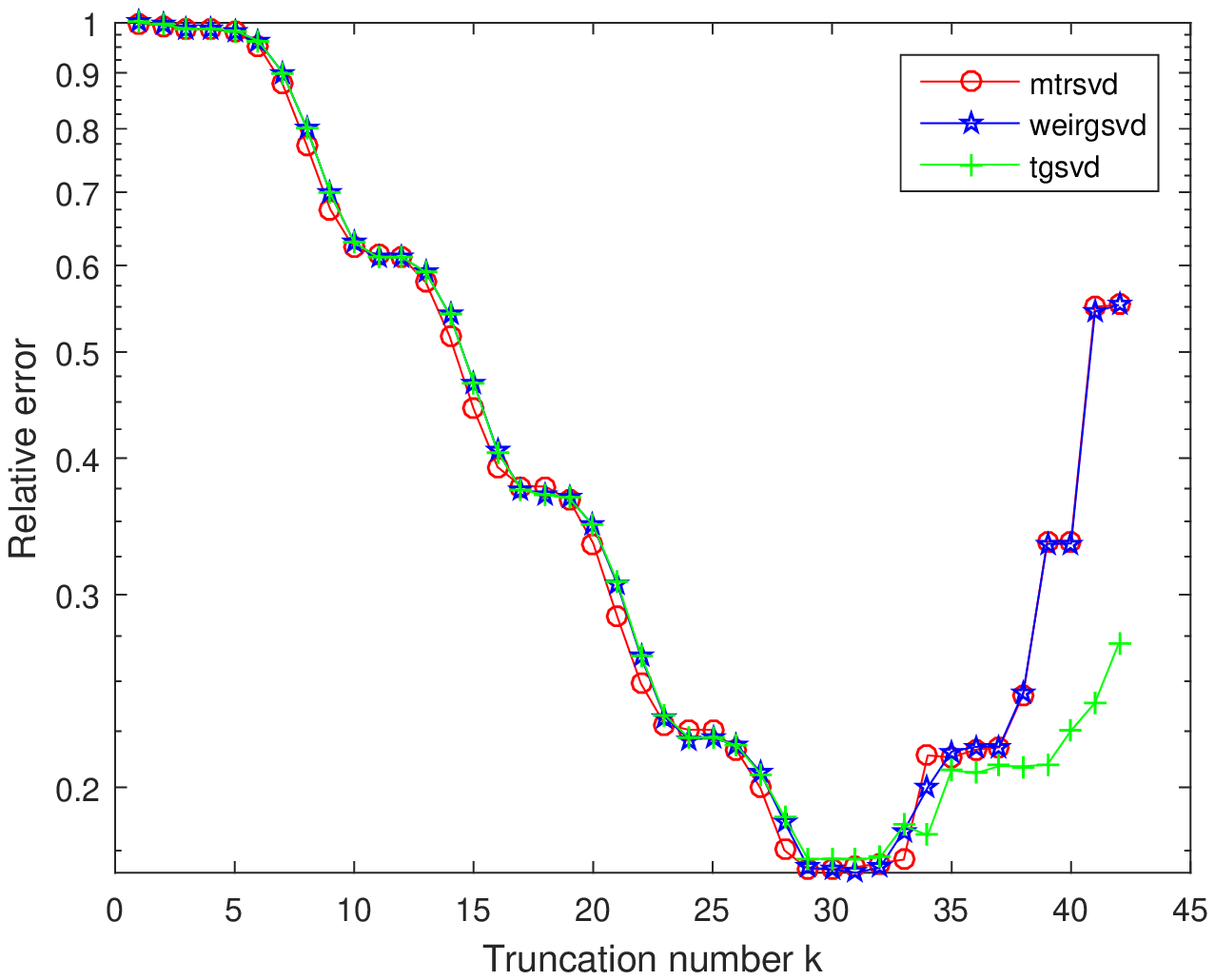}}
  \centerline{(c)}
\end{minipage}
\hfill
\begin{minipage}{0.48\linewidth}
  \centerline{\includegraphics[width=6.0cm,height=4cm]{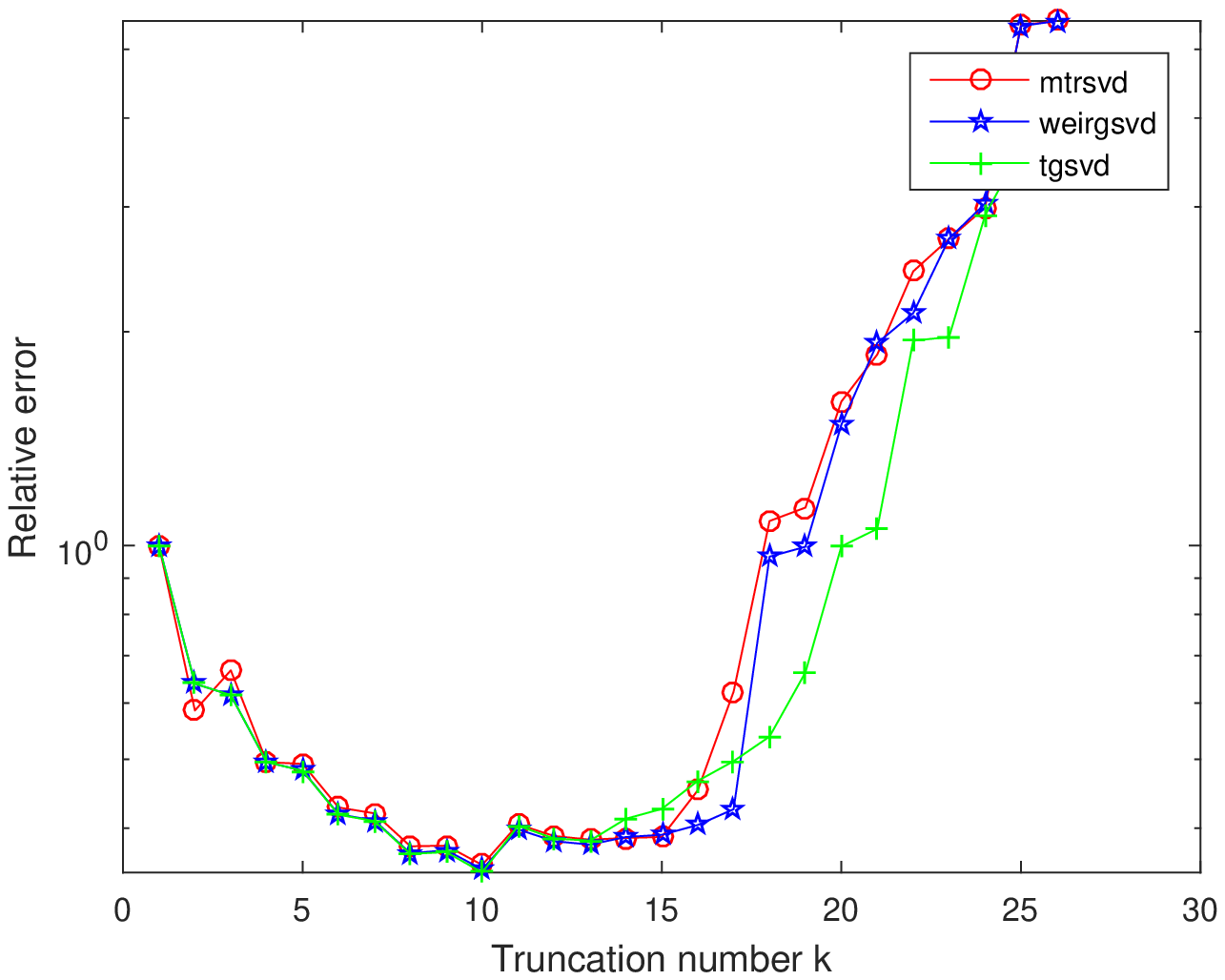}}
  \centerline{(d)}
\end{minipage}
\caption{The semi-convergence processes of  {\sf mtrsvd}, {\sf tgsvd}
and {\sf weirgsvd} for the four test problems with $L=L_3$,
$\varepsilon=10^{-3}$ and $m=n=1,024$: (a) {\sf shaw}; (b) {\sf gravity};
(c) {\sf heat}; (d) {\sf deriv2}.}
\label{fig1}
\end{figure}

Figure \ref{fig1} depicts the convergence processes of {\sf mtrsvd},
{\sf tgsvd} and {\sf weirgsvd} as $k$ increases
for the four test problems with $L=L_3$, $\varepsilon=10^{-3}$ and $m=n=1,024$.
We can see that the three algorithms have very similar convergence processes
and the relative errors of regularized solutions obtained by {\sf mtrsvd} are
almost identical to those by {\sf tgsvd} and {\sf weirgsvd} as $k$ increases
until the occurrence of semi-convergence.
For the other problems, we have observed similar phenomena. These
indicate that the three algorithms
have the same or highly competitive regularizing effects.

\begin{figure}
\begin{minipage}{0.48\linewidth}
  \centerline{\includegraphics[width=6.0cm,height=4cm]{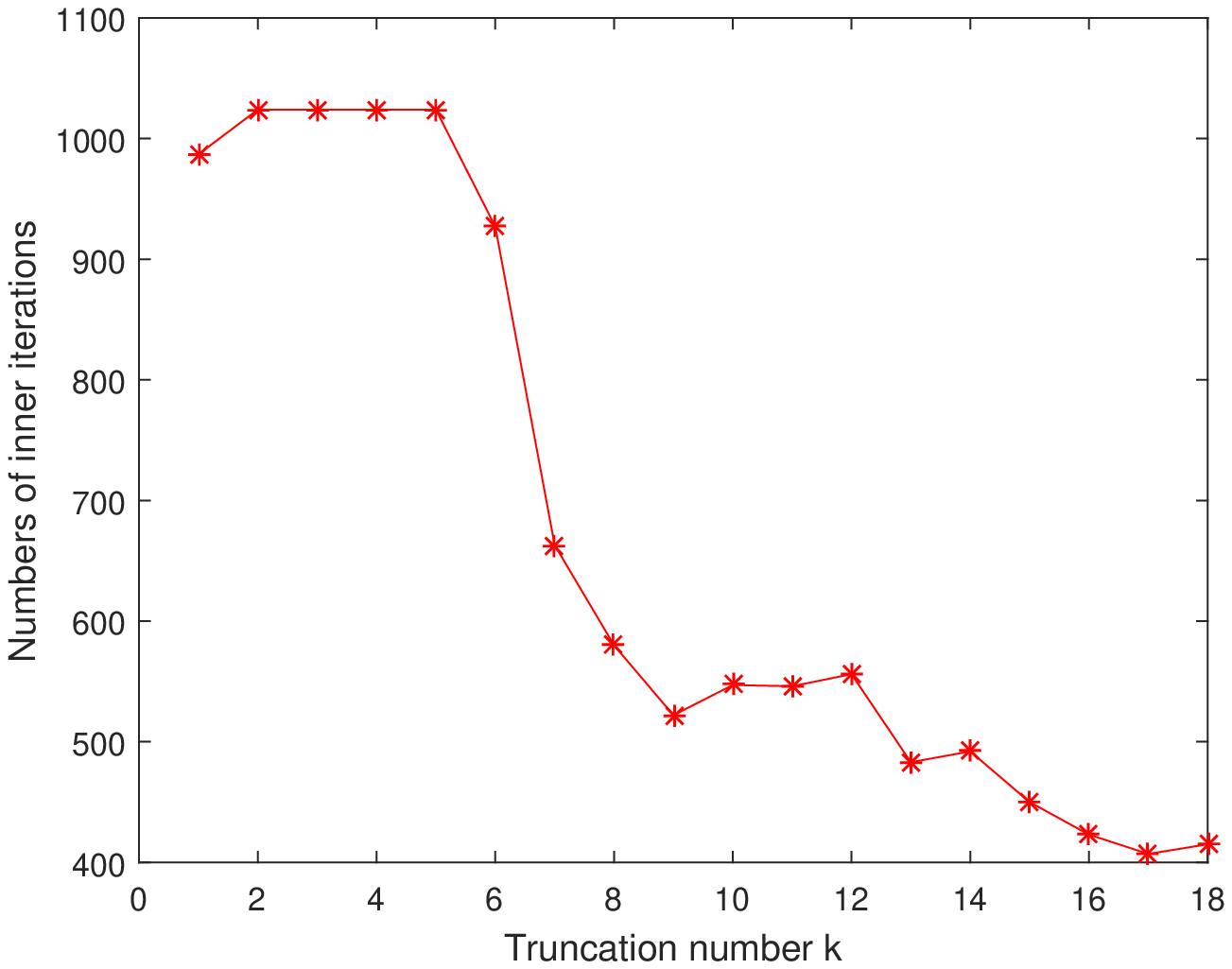}}
  \centerline{(a)}
\end{minipage}
\hfill
\begin{minipage}{0.48\linewidth}
  \centerline{\includegraphics[width=6.0cm,height=4cm]{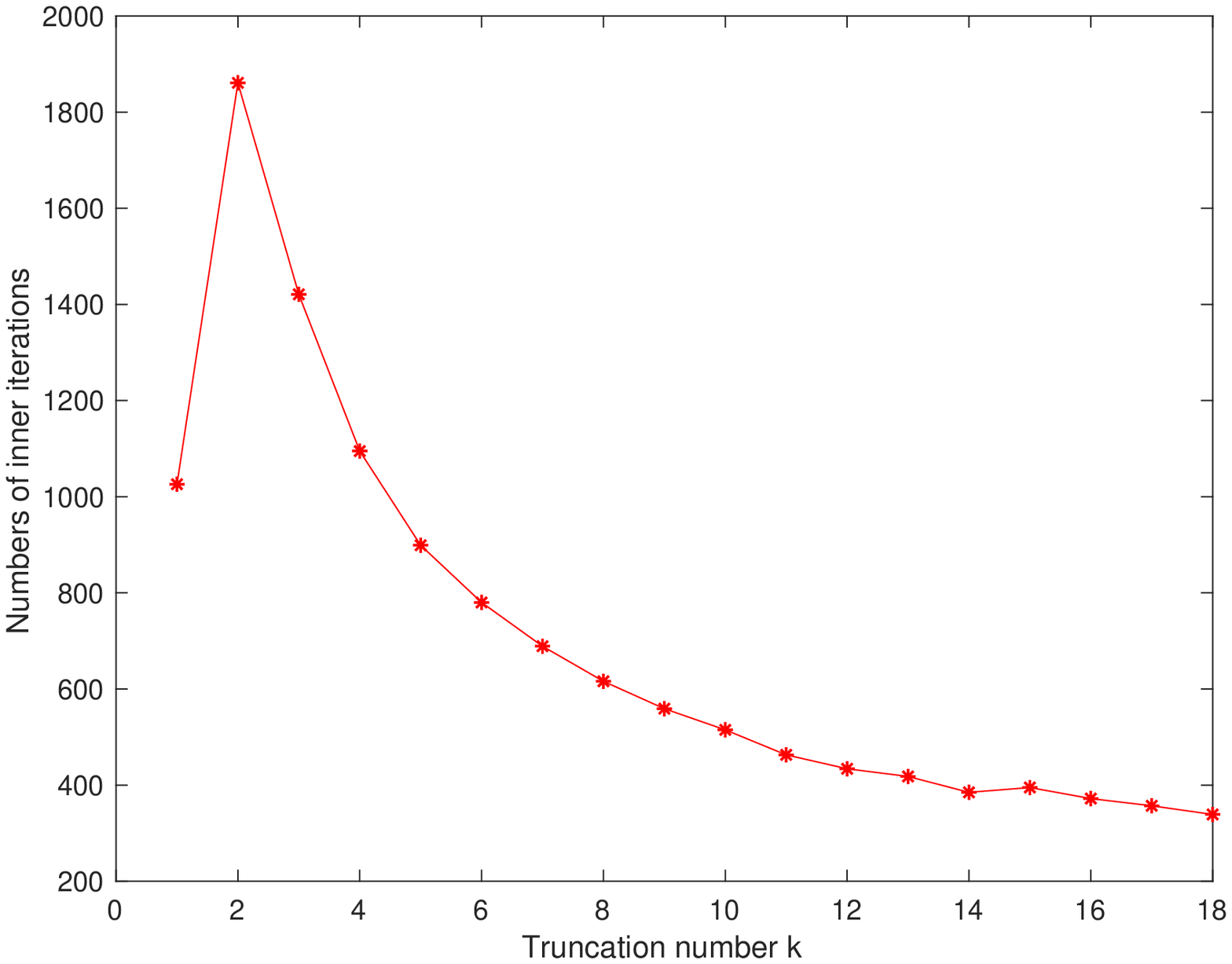}}
  \centerline{(b)}
\end{minipage}

\vfill
\begin{minipage}{0.48\linewidth}
  \centerline{\includegraphics[width=6.0cm,height=4cm]{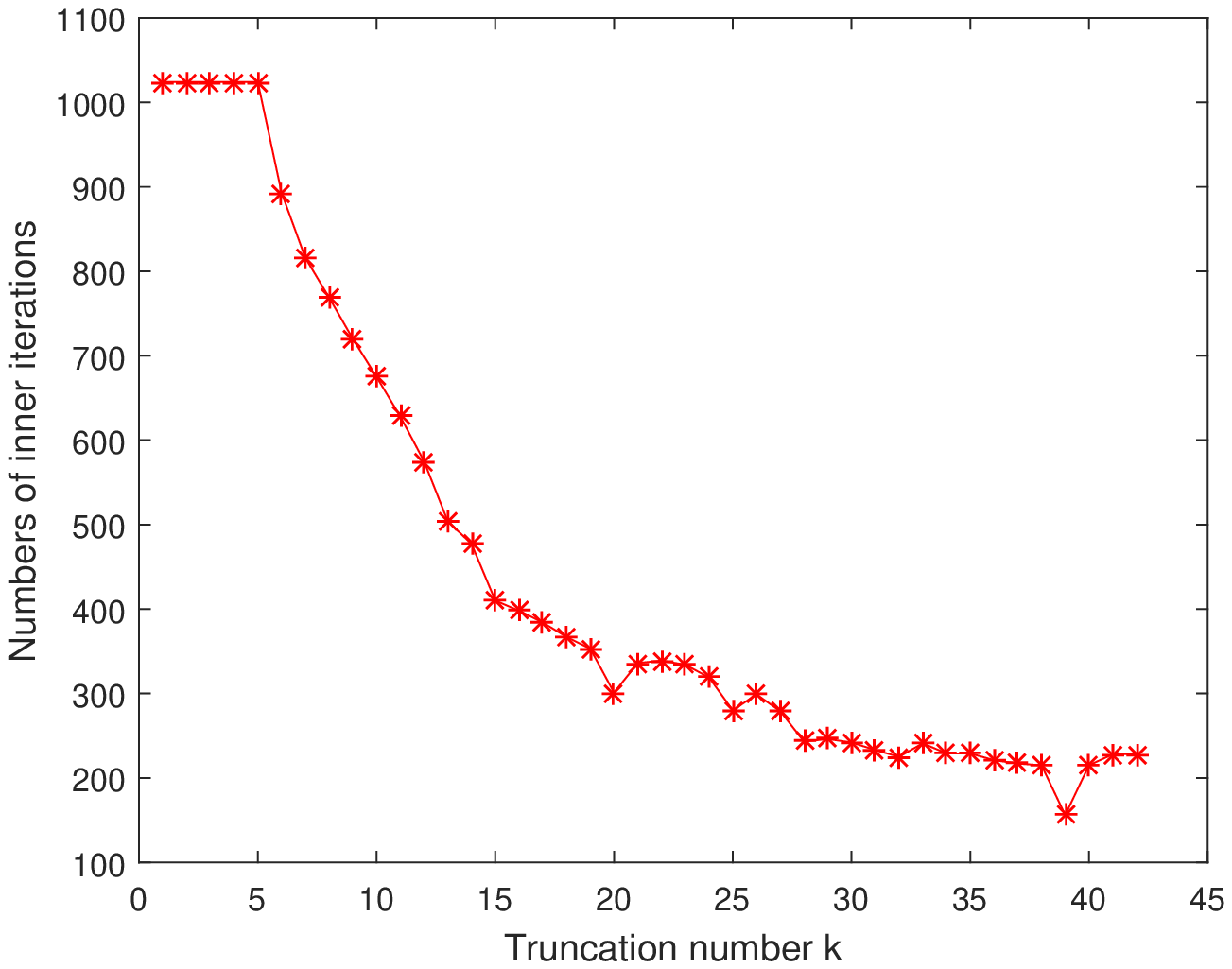}}
  \centerline{(c)}
\end{minipage}
\hfill
\begin{minipage}{0.48\linewidth}
  \centerline{\includegraphics[width=6.0cm,height=4cm]{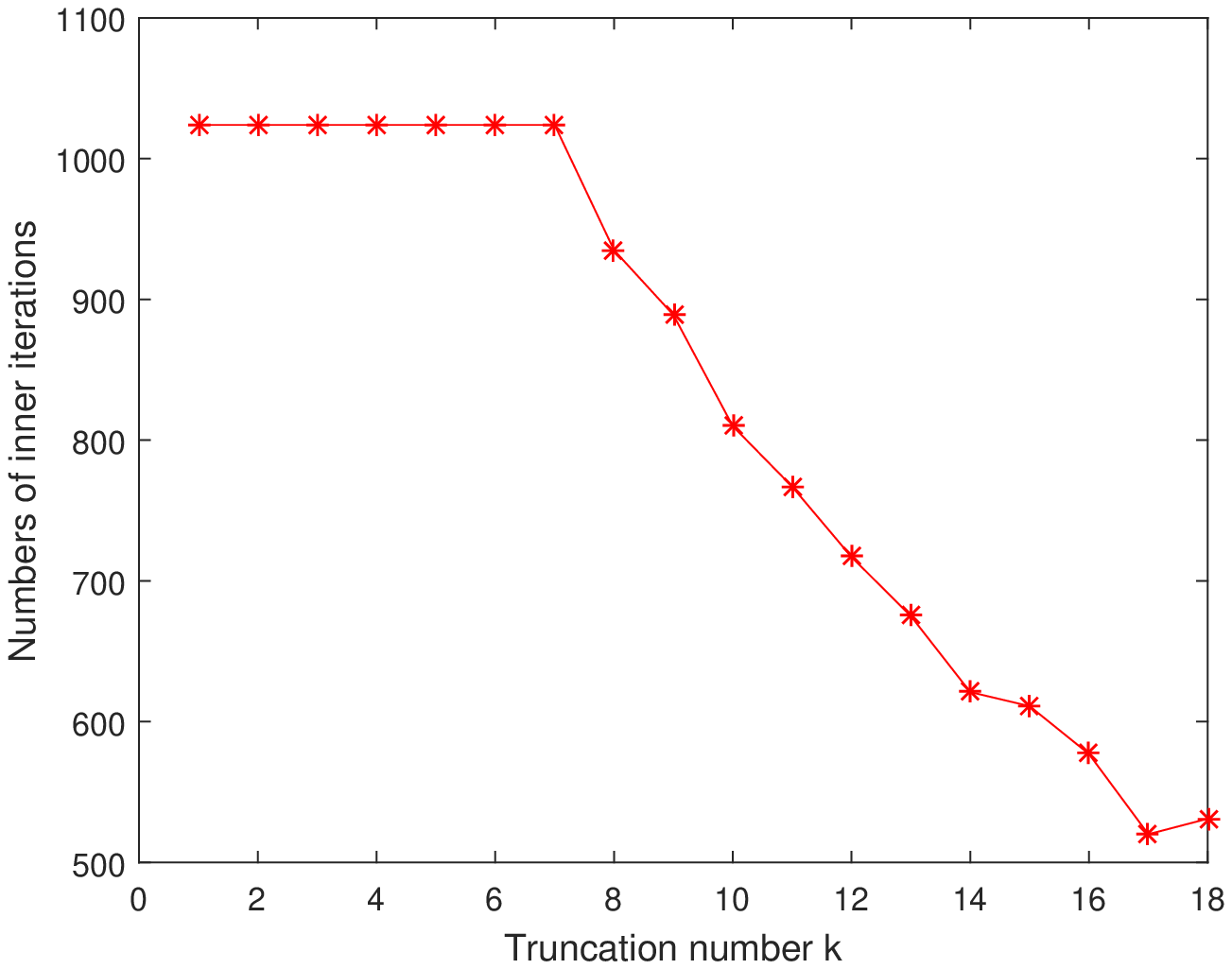}}
  \centerline{(d)}
\end{minipage}
\caption{ The numbers of inner iterations versus $k$ for
Algorithm \ref{alg4} ({\sf mtrsvd}) with $L=L_3$,
$\varepsilon=10^{-3}$ and $m=n=1,024$: (a) {\sf shaw}; (b) {\sf gravity};
(c) {\sf heat}; (d) {\sf deriv2}.}
\label{fig2}
\end{figure}

Figure \ref{fig2} depicts the number of inner iterations used by LSQR
versus the parameter $k$. We clearly observe that the
number of inner iterations exhibits a
considerable decreasing tendency as $k$ increases for the chosen
test problems with $L=L_3$, $\varepsilon=10^{-3}$ and $m=n=1,024$.
LSQR becomes substantially more efficient with $k$ increasing.
For $L=L_1$, we have similar findings.
A distinction is that, for each problem, LSQR uses fewer inner iterations
to converge for $L_1$ than for $L_3$.


\subsubsection{The two dimensional case}

In this subsection, we test the problem
{\sf seismicwavetomo} which is from \cite{hansen12} and creates a
two dimensional seismic tomography.
We use the code of \cite{hansen12} to
generate an $ps \times N^2$ coefficient matrix $A$,
the true solution $x_{true}$
and noise-free right-hand side $b_{true}$.
We take $N=32$ and $100$ with default
$s=N$ and $p=2N$, respectively, that is, we generate
$A\in \mathbb{R}^{2,048\times 1,024}$ and $A\in
\mathbb{R}^{20,000\times 10,000}$.

\begin{table}[h]\small
  \centering
  \caption{The relative errors for
  {\sf seismicwavetomo}.}
    \label{tab6}
    \centerline{$m=2,048$ and $n=1,024$}
    \begin{minipage}[t]{1\textwidth}
     \begin{tabular*}{\linewidth}{lp{0.1cm}p{1.475cm}p{1.47cm}p{1.47cm}p{1.47cm}p{1.47cm}p{1.47cm}p{1.47cm}p{1.47}p{1.47cm}}
     \toprule[0.6pt]
     \multicolumn{2}{c}{}       &\multicolumn{3}{c}{$L=L_1$}  &\multicolumn{3}{c}{$L=L_3$} \\
     \cmidrule(lr){3-5}\cmidrule(lr){6-8}
     $\varepsilon$ &$q$  &{\sf tgsvd}       &{\sf weirgsvd}	    &{\sf mtrsvd}	
     &{\sf tgsvd}         &{\sf weirgsvd}      &{\sf mtrsvd}  \\ \midrule[0.6pt]
     $10^{-2}$     &70 &0.6397(347)	&0.6358(323)	&0.6116(305)  &0.7526(369)	 &0.7498(380)    &0.7181(285)\\
     $10^{-3}$     &42 &0.3117(658)	&0.3083(603)	&0.2982(586)  &0.3777(623)	 &0.3691(633)    &0.3451(604)\\
     \bottomrule[0.6pt]
     \end{tabular*}\\[2pt]
     \end{minipage}

  \centerline{$m=20,000$ and $n= 10,000$}

     \begin{minipage}[t]{1\textwidth}
     \begin{tabular*}{\linewidth}{lp{2.0cm}p{2.0cm}p{2.0cm}p{1.0cm}p{1.0cm}}
     \toprule[0.6pt]
     \multicolumn{1}{c}{}       &\multicolumn{2}{c}{$L=L_1$}  &\multicolumn{2}{c}{$L=L_3$} \\
     \cmidrule(lr){2-3}\cmidrule(lr){4-5}
     $\varepsilon$ &$q$   &{\sf mtrsvd}	        &$q$       &{\sf mtrsvd}  \\ \midrule[0.6pt]
     $10^{-2}$     &141   &0.8691(419)	     &295       &0.9092(345)  \\
     $10^{-3}$     &951   &0.7766(1249)	    &53	       &0.8949(1147)  \\
     \bottomrule[0.6pt]
     \end{tabular*}\\[2pt]
  \end{minipage}

\end{table}

Table \ref{tab6} shows the relative errors of the best regularized solutions
obtained by  {\sf mtrsvd}, {\sf tgsvd} and {\sf weirgsvd} with
$L=L_1, L_3$ and $\varepsilon=10^{-2},10^{-3}$, respectively, where
{\sf tgsvd} and {\sf weirgsvd} are out of memory for $m=20,000,n=10,000$.
Obviously, the relative errors of the best regularized solutions by {\sf mtrsvd}
are at least as accurate as those by {\sf tgsvd}
and {\sf weirgsvd} for the two given $\varepsilon$ and $m=2,048, n=1,024$,
and {\sf mtrsvd} is more practical than {\sf tgsvd} and {\sf weirgsvd}
for large scale problems.

\begin{table}[h]\small
  \centering
  \caption{The relative errors for
  {\sf seismicwavetomo} of $m=20,000, n= 10,000$.}
    \label{addtab1}
    \begin{minipage}[t]{1\textwidth}
     \begin{tabular*}{\linewidth}{lp{2.0cm}p{2.0cm}p{2.0cm}p{1.0cm}}
     \toprule[0.6pt]
     \multicolumn{0}{c}{} &\multicolumn{2}{c}{$L=L_1,\varepsilon = 10^{-2}$}
     &\multicolumn{2}{c}{$L=L_3,\varepsilon = 10^{-3} $} \\
  \cmidrule(lr){2-3}\cmidrule(lr){4-5}
           &$q$       &{\sf mtrsvd}           &$q$       &{\sf mtrsvd}\\ \midrule[0.6pt]
           &75        &0.9489(405)	    &10        &0.7907(1110)  \\
           &173       &0.9483(407)  	&587       &0.7802(1213) \\
           &176       &0.9475(424)      &951       &0.7766(1249)\\
     \bottomrule[0.6pt]
     \end{tabular*}\\[2pt]
  \end{minipage}
\end{table}

We next investigate how {\sf mtrsvd} behaves as the oversampling parameter
 $q$ varies for this problem with $m=20,000$ and $n=10,000$.
 Table \ref{addtab1} shows the relative errors of the
 best regularized solutions
obtained by {\sf mtrsvd} for varying $q$ with $L=L_1, L_3$ and
$\varepsilon=10^{-2},10^{-3}$, respectively. As we can see,
the relative errors of the best regularized solutions by {\sf mtrsvd}
for {\sf seismicwavetomo} decrease a little bit with $q$ increasing.
This confirms our theory that bigger $q$ should generally generate more
accurate rank-$k$ approximation to $A$, so that the regularized
solutions could be more accurate.
%

%

\begin{figure}
\begin{minipage}{0.48\linewidth}
  \centerline{\includegraphics[width=6.0cm,height=4cm]{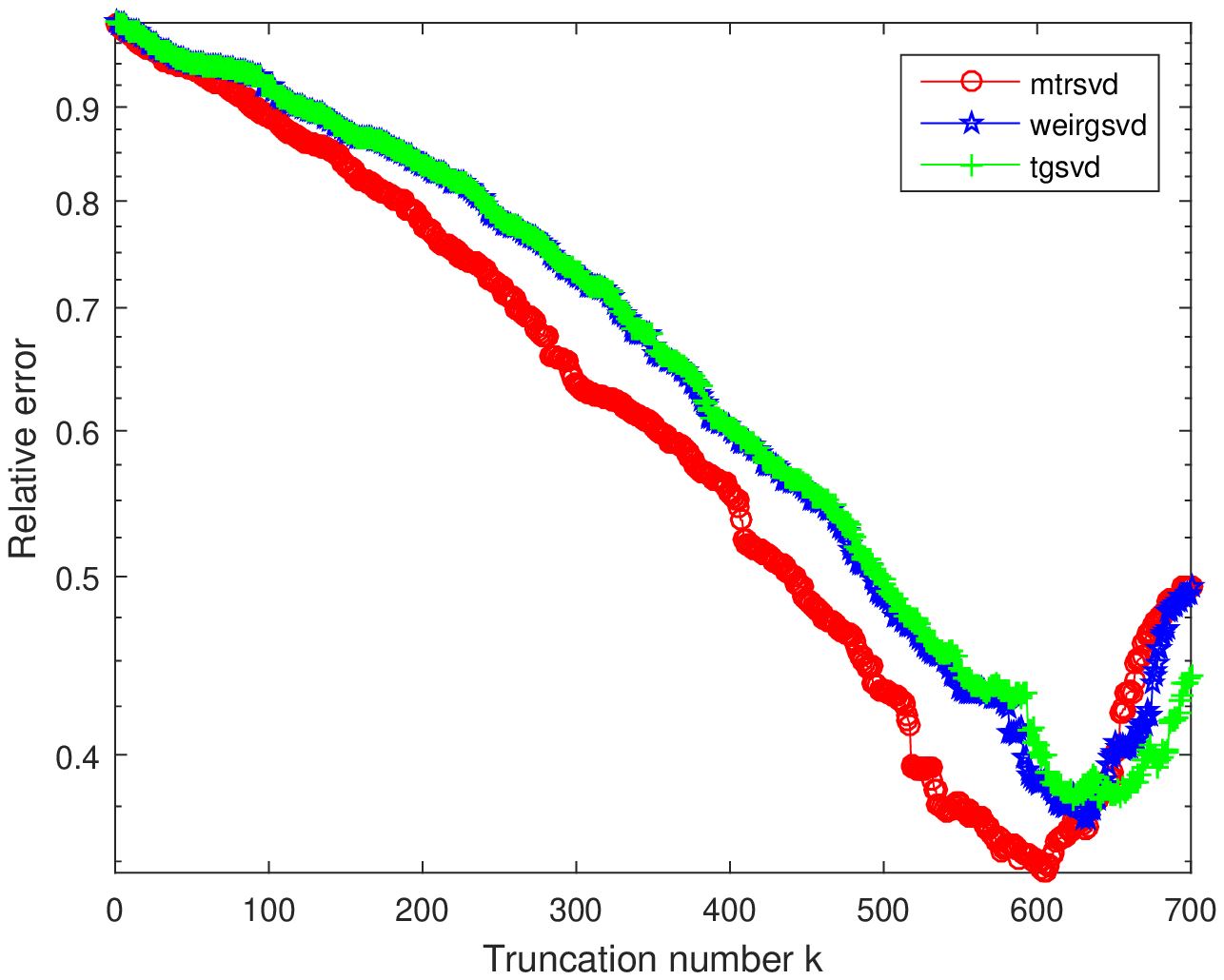}}
\end{minipage}
\hfill
\begin{minipage}{0.48\linewidth}
  \centerline{\includegraphics[width=6.0cm,height=4cm]{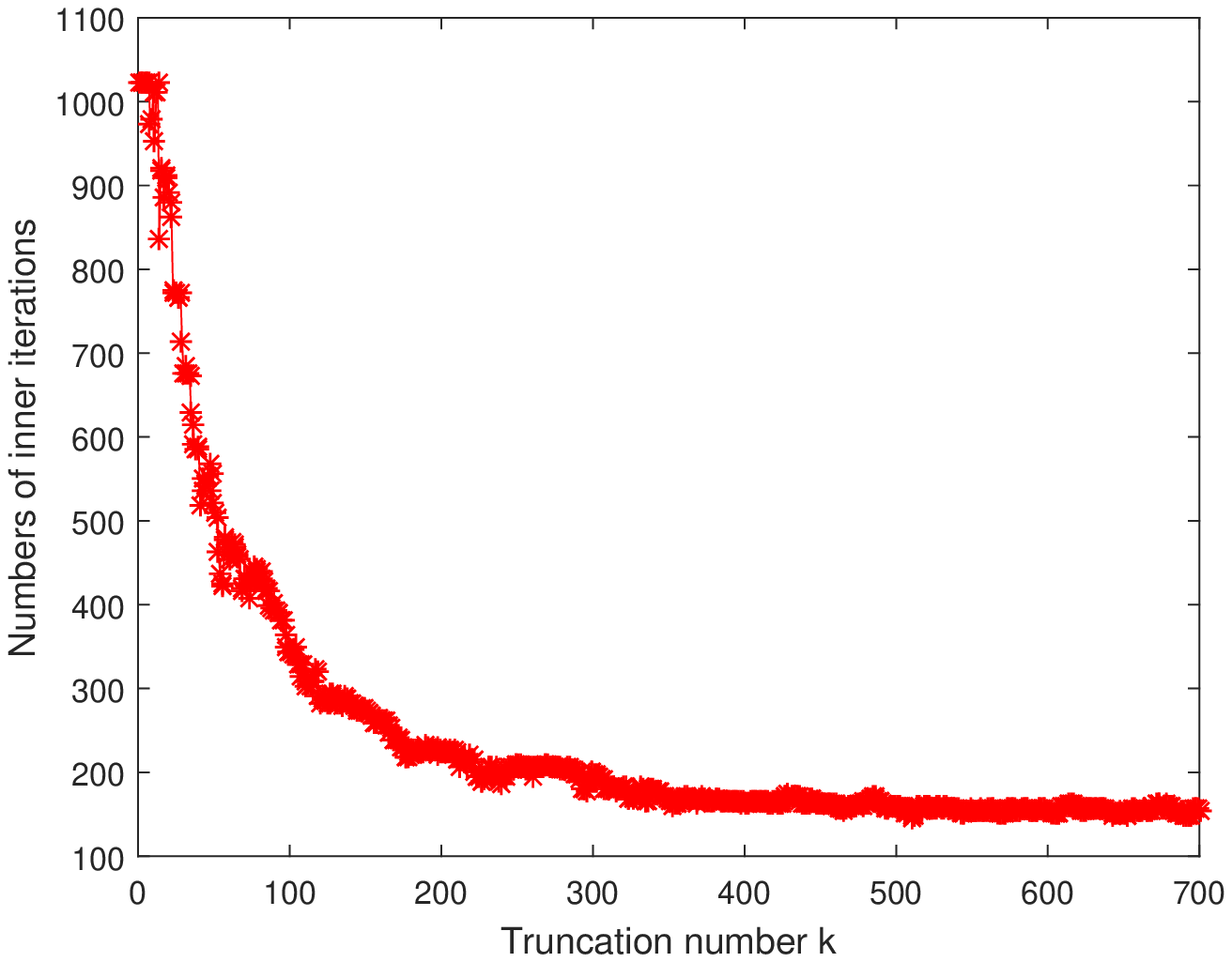}}
\end{minipage}
\caption{The relative errors of {\sf tgsvd}, {\sf weirgsvd} and {\sf mtrsvd},
and inner iterations versus $k$
of Algorithm~\ref{alg4} ({\sf mtrsvd}) for the problem
{\sf seismicwavetomo} of $m=2,048, n=1,024$
with $L=L_3$ and $\varepsilon=10^{-3}$.}
\label{addfig2}
\end{figure}

Figures \ref{addfig2} draws the convergence processes of {\sf mtrsvd}
{\sf tgsvd} and {\sf weirgsvd} for $m=20,48, n=10,24$ and the
inner iterations versus the parameter $k$ with $\varepsilon=10^{-3}$ and
$L=L_3$. We can see that the best regularized solution by {\sf mtrsvd}
is more accurate than the counterparts by {\sf tgsvd}
and {\sf weirgsvd} and LSQR uses substantially fewer
iterations as $k$ increases.
Compared with the results on the one dimensional problems, however, we observe
a remarkable difference that the optimal regularization parameter $k_0$ now
becomes much bigger. The reason is that for this problem, as we have
numerically justified by plotting the discrete Picard condition, the Fourier
coefficients $| u_i^T b_{true}|$ do not decay considerably faster than the
generalized singular values $\sigma_i$ of $\{A,L,\}$, where the $u_i$ are the
first $\min\{p,n\}$ left singular vectors of $\{A,L\}$. Recall that the GSVD of
$\{A,L\}$ is $A=UCZ^{-1}$ and $L=VSZ^{-1}$, where $C \in \mathbb{R}^{m \times n}$
and $S \in \mathbb{R}^{p \times n}$ are diagonal matrices with
the diagonal entries $c_i$ and $s_i$, respectively,
$C^TC+S^TS=I$, $\sigma_i=c_i/s_i$, $U\in \mathbb{R}^{m \times m}$ and
$V \in \mathbb{R}^{p \times p}$ are orthogonal, and
the columns $u_i$ of $U$ are called the left singular vectors.
This means that a good regularized solution
must include many dominant GSVD components of $\{A,L\}$.

\subsection{The $m\leq n$ case}

We now test Algorithm~\ref{alg5} ({\sf mtrsvd}), {\sf tgsvd}, {\sf weirgsvd}
and {\sf xiangrgsvd} on the test problems in Table \ref{tab1}.
In Table~\ref{tab7}, we display the relative errors of the best regularized
solutions obtained by {\sf mtrsvd}, {\sf tgsvd}, {\sf weirgsvd} and {\sf xiangrgsvd}
with $L = L_1$ and $\varepsilon =10^{-2}, 10^{-3}$, respectively.

\begin{table}[htp]
     \centering
     \caption{The comparison of Algorithm \ref{alg5} ({\sf mtrsvd}) and the others
     with $L=L_1$ and
    $\varepsilon=10^{-2}, \ 10^{-3}$.}
     \label{tab7}
     \centerline{$\varepsilon=10^{-2}$}
     \begin{minipage}[t]{1\textwidth}
     \begin{tabular*}{\linewidth}{lp{0.1cm}p{1.4cm}p{1.4cm}p{1.4cm}p{1.4cm}p{1.4cm}p{1.4cm}p{1.4cm}p{1.4cm}p{1.4cm}}
     \toprule[0.6pt]
    \multicolumn{2}{c}{}   &\multicolumn{4}{c}{$m=n=1,024$} &\multicolumn{2}{c}{$m=n=10,240$} \\
     \cmidrule(lr){3-6}\cmidrule(lr){7-8}
     	            &$q$   &{\sf tgsvd} &{\sf xiangrgsvd}	    &{\sf weirgsvd}		&{\sf mtrsvd}
       	     &{\sf xiangrgsvd}   &{\sf mtrsvd}\\  \midrule[0.6pt]
     {\sf shaw}     &11   &0.2099(6)     &0.2099(7)	  &0.2099(7)	 &0.2097(7)      &0.1666(8)	  &0.1669(8)\\
     {\sf gravity}  &9   &0.3004(8)	     &0.2993(9)	  &0.2993(9)	 &0.2993(9)     &0.2743(10)  &0.2785(10)\\
     {\sf heat}     &7   &0.2228(27)	&0.3561(23)	  &0.3561(23)	 &0.2488(24)     &0.3535(23)  &0.2369(25)\\
     {\sf deriv2}   &12   &0.4329(5)	&3.1031(1)	  &3.1031(1)	 &0.4455(6)      &3.1025(1)	  &0.4343(8)\\
     \bottomrule[0.6pt]
     \end{tabular*}\\[2pt]
  \end{minipage}
     \centerline{$\varepsilon=10^{-3}$}
  \begin{minipage}[t]{1\textwidth}
  \begin{tabular*}{\linewidth}{lp{0.1cm}p{1.4cm}p{1.4cm}p{1.4cm}p{1.4cm}p{1.4cm}p{1.5cm}p{1.4cm}p{1.4cm}p{1.4cm}}
  \toprule[0.6pt]
  \multicolumn{2}{c}{}   &\multicolumn{4}{c}{$m=n=1,024$} &\multicolumn{2}{c}{$m=n=10,240$} \\
     \cmidrule(lr){3-6}\cmidrule(lr){7-8}
     	            &$q$   &{\sf tgsvd}        &{\sf xiangrgsvd}	    &{\sf weirgsvd}
     &{\sf mtrsvd}  	     &{\sf xiangrgsvd}   &{\sf mtrsvd}\\  \midrule[0.6pt]
     {\sf shaw}     &6   &0.1946(6)     &0.1967(7)	  &0.1967(7)	 &0.1942(7)      &0.1353(9)	  &0.1311(9)\\
     {\sf gravity}  &6   &0.2556(12)	&0.2382(12)	  &0.2382(12)	 &0.2577(11)     &0.2443(12)  &0.2223(12)\\
     {\sf heat}     &8   &0.1543(30)	&0.1714(30)	  &0.1714(30)	 &0.1564(29)     &0.1801(32)  &0.1523(32)\\
     {\sf deriv2}   &6   &0.3342(12)	&3.0987(1)	  &3.0987(1)	 &0.3790(8)      &3.0985(1)	  &0.3815(8)\\
  \bottomrule[0.6pt]
  \end{tabular*}\\[2pt]
  \end{minipage}
\end{table}

The results indicate that for $m=n=1,024$ {\sf mtrsvd}
computes the best regularized solution with very similar accuracy to those by
{\sf tgsvd} and {\sf weirgsvd}, {\sf xiangrgsvd} for severely and
moderately ill-posed problems, but the solution accuracy by {\sf mtrsvd} is
much higher than that by {\sf xiangrgsvd} and {\sf weirgsvd}
for the mildly ill-posed problem {\sf deriv2}. Actually, the
best regularized solutions by {\sf xiangrgsvd} and {\sf weirgsvd} have no accuracy
since their relative errors are over 300\%! As is expected,
whenever an algorithm has regularizing
effects and can compute a regularized solution with some accuracy,
the smaller $\varepsilon$ is, the bigger $k_0$ is and the more
accurate regularized solution is, except for {\sf shaw} of
$m=1,024$ and $n=1,024$ where
the $k_0$ are the same for each algorithm with $\varepsilon=10^{-2},10^{-3}$.

Mathematically, {\sf weirgsvd} is the same as {\sf xiangrgsvd}.
Table \ref{tab7} confirms that these two algorithms compute the
same regularized solutions for $m=n=1,024$.
For this reason, we only report the
results obtained by {\sf xiangrgsvd} for $m=n=10,240$. Still, {\sf xiangrgsvd} fails
to solve {\sf deriv2} and the relative errors of the best regularized
solutions are over 300\%, but {\sf mtrsvd} is successful to obtain
good regularized solutions.

\begin{table}[htp]
  \centering
  \caption{The comparison of Algorithm \ref{alg5} ({\sf mtrsvd})
  and the others with $L=L_3$ and
   $\varepsilon=10^{-2},\ 10^{-3}$.}
  \label{tab9}
  \centerline{$\varepsilon=10^{-2}$}
  \begin{minipage}[t]{1\textwidth}
  \begin{tabular*}{\linewidth}{lp{0.1cm}p{1.4cm}p{1.4cm}p{1.4cm}p{1.4cm}p{1.4cm}p{1.5cm}p{1.4cm}p{1.4cm}p{1.4cm}}
  \toprule[0.6pt]
     \multicolumn{2}{c}{}   &\multicolumn{4}{c}{$m=n=1,024$} &\multicolumn{2}{c}{$m=n=10,240$} \\
     \cmidrule(lr){3-6}\cmidrule(lr){7-8}
     	            &$q$   &{\sf tgsvd}  &{\sf xiangrgsvd}	    &{\sf weirgsvd}		
     &{\sf mtrsvd}  	     &{\sf xiangrgsvd}   &{\sf mtrsvd}\\ \midrule[0.6pt]
     {\sf shaw}     &9   &0.2000(7)	   &0.2000(7)	    &0.2000(7)	    &0.1993(7)       &0.1475(9)	   &0.1371(9)\\
     {\sf gravity}  &10   &0.3287(8)   &0.3265(8)	    &0.3265(8)	    &0.3280(8)	     &0.2315(12)   &0.1787(12)\\
     {\sf heat}     &9  &0.3065(24)   &0.3721(20)	    &0.3721(20)	    &0.3268(21)	     &0.1806(35)   &0.1402(35)\\
     {\sf deriv2}   &12   &0.4481(6)   &3.1009(1)	    &3.1009(1)	    &0.4905(5)       &3.0993(1)	   &0.3758(9)\\
    \bottomrule[0.6pt]
    \end{tabular*}\\[2pt]
  \end{minipage}
  \centerline{$\varepsilon=10^{-3}$}
  \begin{minipage}[t]{1\textwidth}
  \begin{tabular*}{\linewidth}{lp{0.1cm}p{1.4cm}p{1.4cm}p{1.4cm}p{1.4cm}p{1.4cm}p{1.5cm}p{1.4cm}p{1.4cm}p{1.4cm}}
  \toprule[0.6pt]
       \multicolumn{2}{c}{}   &\multicolumn{4}{c}{$m=n=1,024$} &\multicolumn{2}{c}{$m=n=10,240$} \\
     \cmidrule(lr){3-6}\cmidrule(lr){7-8}
     	            &$q$   &{\sf tgsvd}        &{\sf xiangrgsvd}	    &{\sf weirgsvd}		
     &{\sf mtrsvd}  	 &{\sf xiangrgsvd}    &{\sf mtrsvd}\\ \midrule[0.6pt]
     {\sf shaw}     &6   &0.1659(8)	   &0.1662(8)	    &0.1662(8)	    &0.1659(8)       &0.2010(7)	   &0.1998(7)\\
     {\sf gravity}  &8   &0.2686(10)   &0.2655(10)	    &0.2655(10)	    &0.2668(10)	     &0.2891(11)   &0.2926(11)\\
     {\sf heat}     &9  &0.1689(31)   &0.2851(30)	    &0.2851(30)	    &0.1825(28)	     &0.4519(21)   &0.3371(22)\\
     {\sf deriv2}   &8   &0.3374(12)   &3.0999(1)	    &3.0999(1)	    &0.3891(8)       &3.0984(1)	   &0.3857(10)\\
    \bottomrule[0.6pt]
    \end{tabular*}\\[2pt]
  \end{minipage}
\end{table}

In Table \ref{tab9}, we display the relative errors of best regularized
solutions by all the algorithms with $L = L_3$ and $\varepsilon =10^{-2},
10^{-3}$, respectively. Clearly, we can observe very similar phenomena to
those in Table \ref{tab7}.

\begin{figure}
\begin{minipage}{0.48\linewidth}
  \centerline{\includegraphics[width=6.0cm,height=4cm]{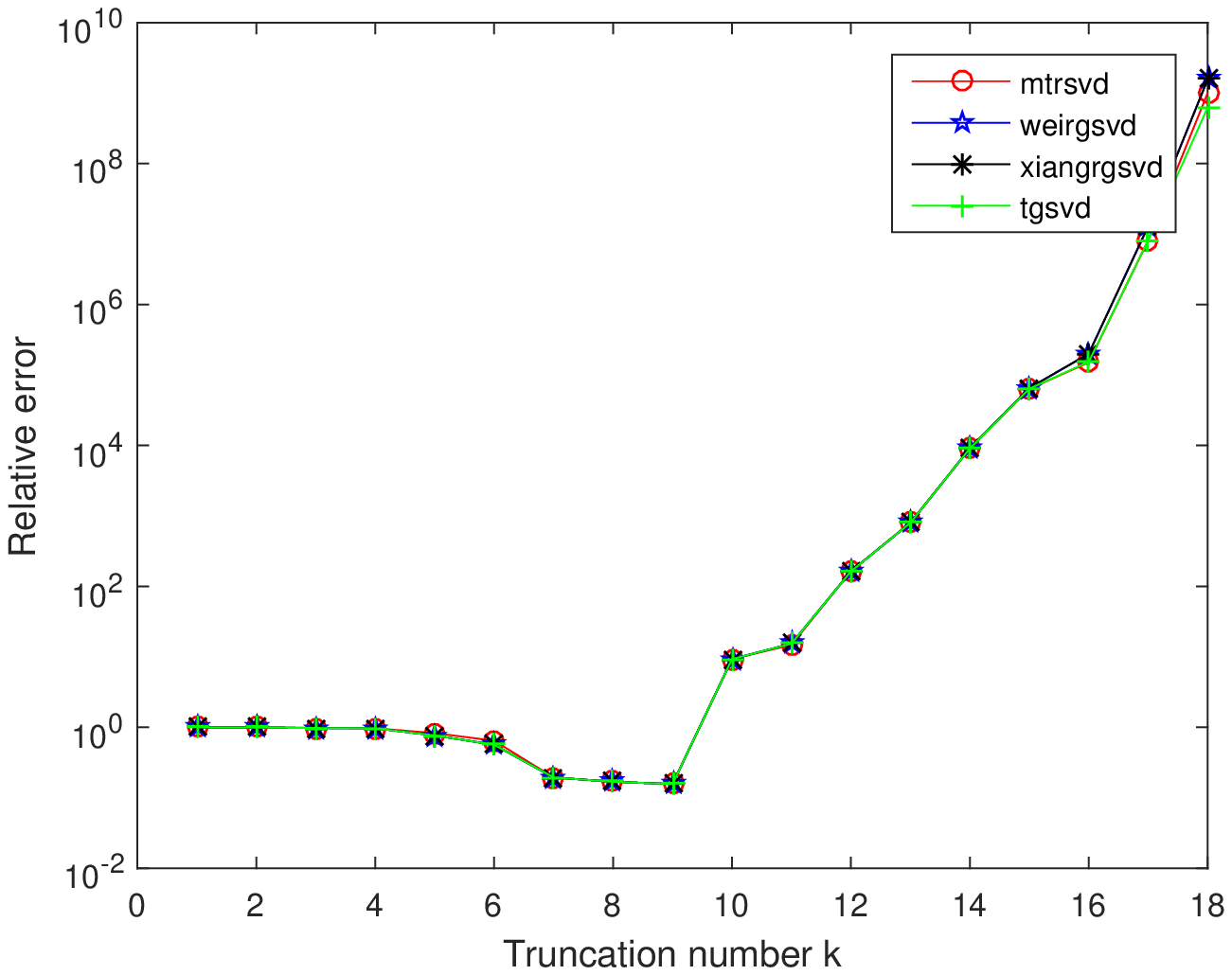}}
  \centerline{(a)}
\end{minipage}
\hfill
\begin{minipage}{0.48\linewidth}
  \centerline{\includegraphics[width=6.0cm,height=4cm]{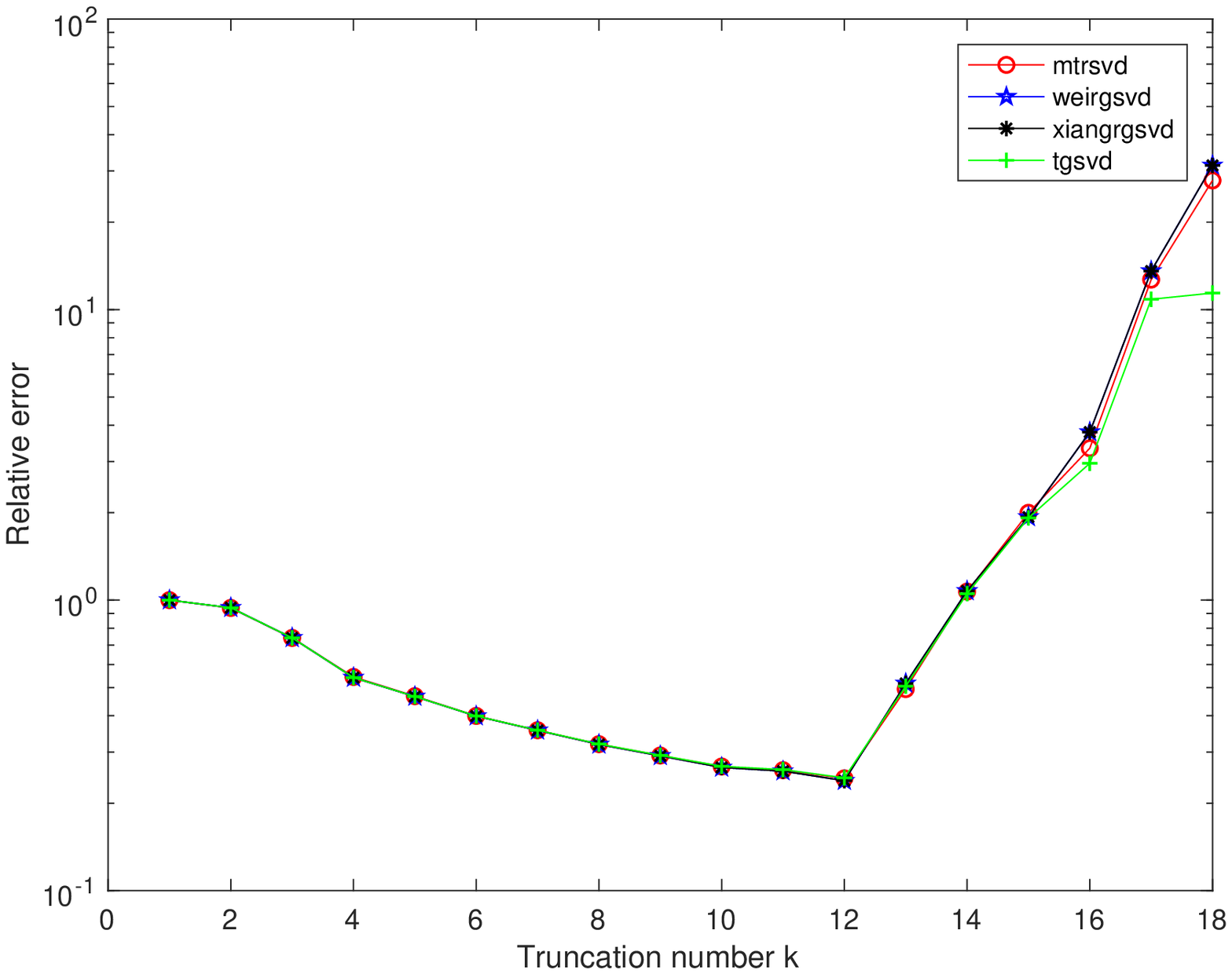}}
  \centerline{(b)}
\end{minipage}

\vfill
\begin{minipage}{0.48\linewidth}
  \centerline{\includegraphics[width=6.0cm,height=4cm]{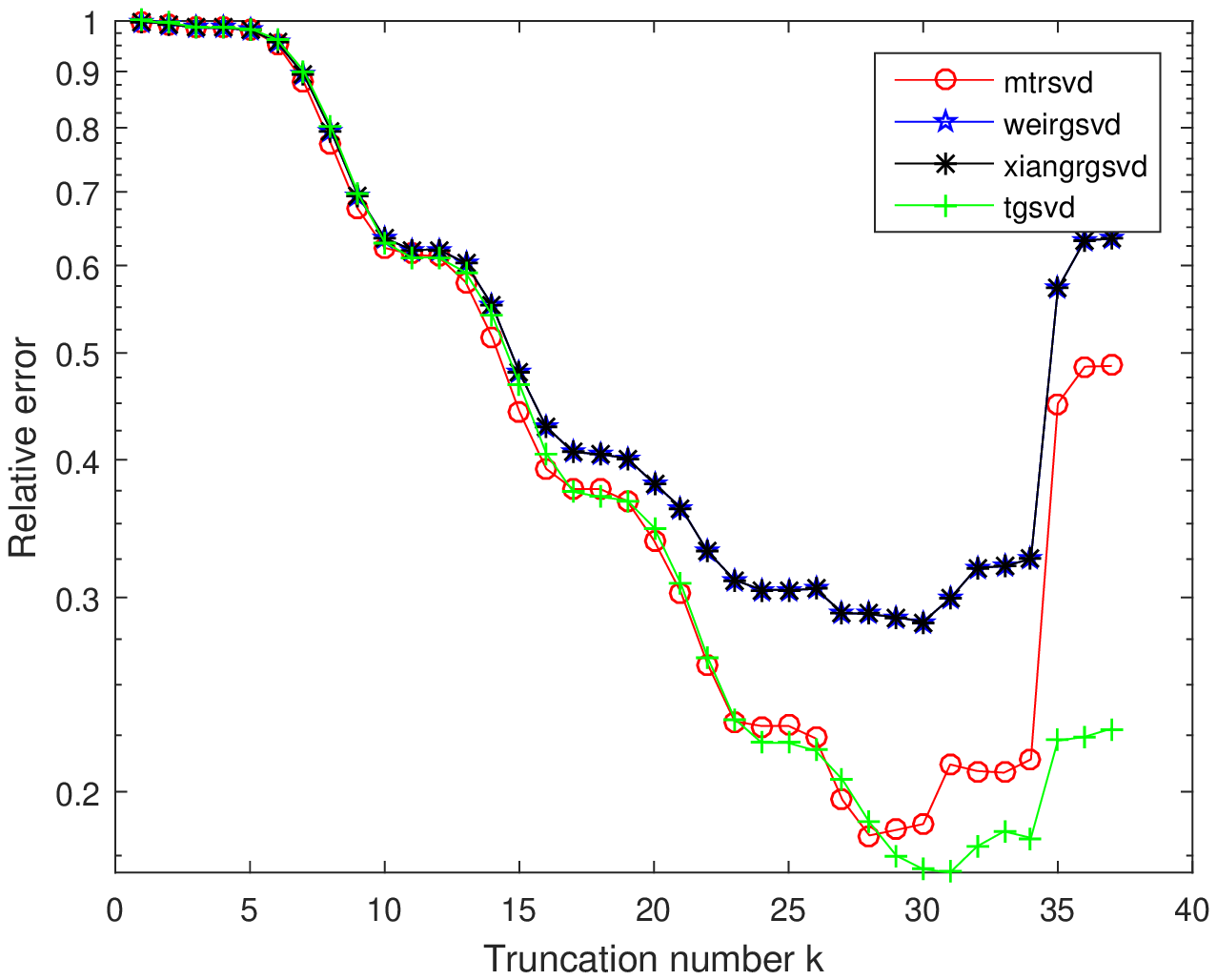}}
  \centerline{(c)}
\end{minipage}
\hfill
\begin{minipage}{0.48\linewidth}
  \centerline{\includegraphics[width=6.0cm,height=4cm]{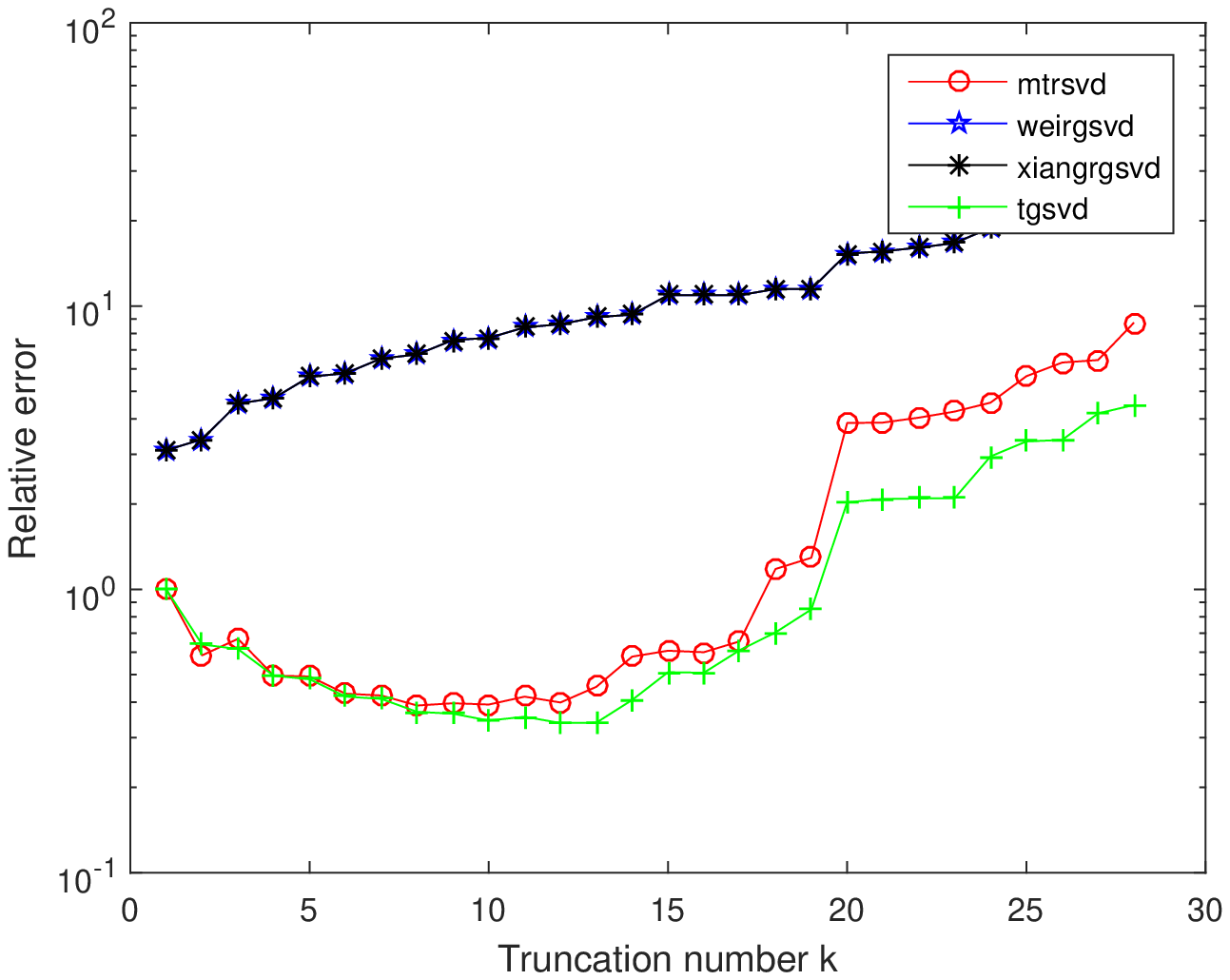}}
  \centerline{(d)}
\end{minipage}
\caption{ The relative errors of Algorithm \ref{alg5} ({\sf mtrsvd}) with $L=L_3$ and
$\varepsilon=10^{-3}$ and $m=n=1,024$: (a) {\sf shaw};
(b) {\sf gravity}; (c) {\sf heat}; (d) {\sf deriv2}.}
\label{fig5}
\end{figure}

\begin{figure}
\begin{minipage}{0.48\linewidth}
  \centerline{\includegraphics[width=6.0cm,height=4cm]{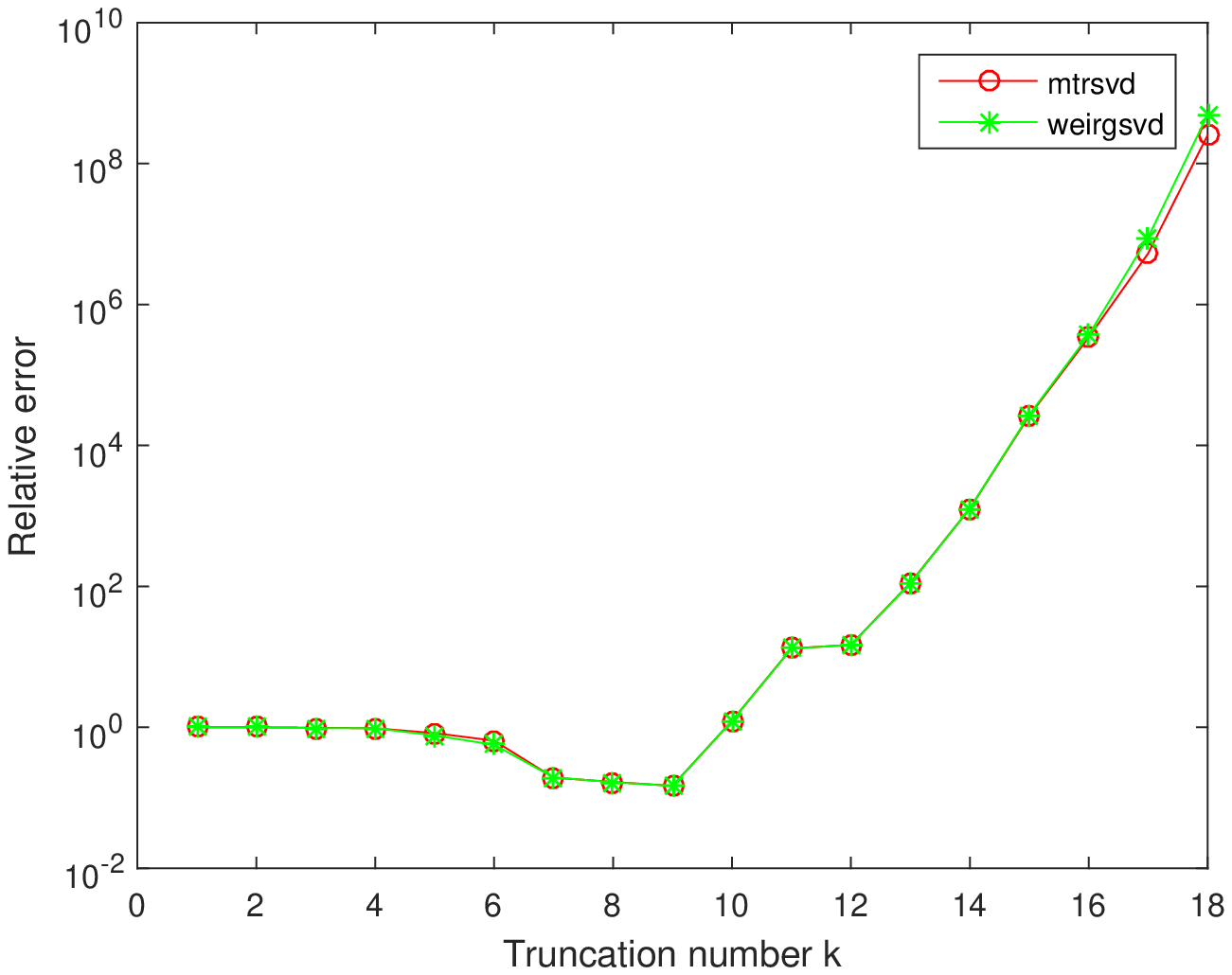}}
  \centerline{(a)}
\end{minipage}
\hfill
\begin{minipage}{0.48\linewidth}
  \centerline{\includegraphics[width=6.0cm,height=4cm]{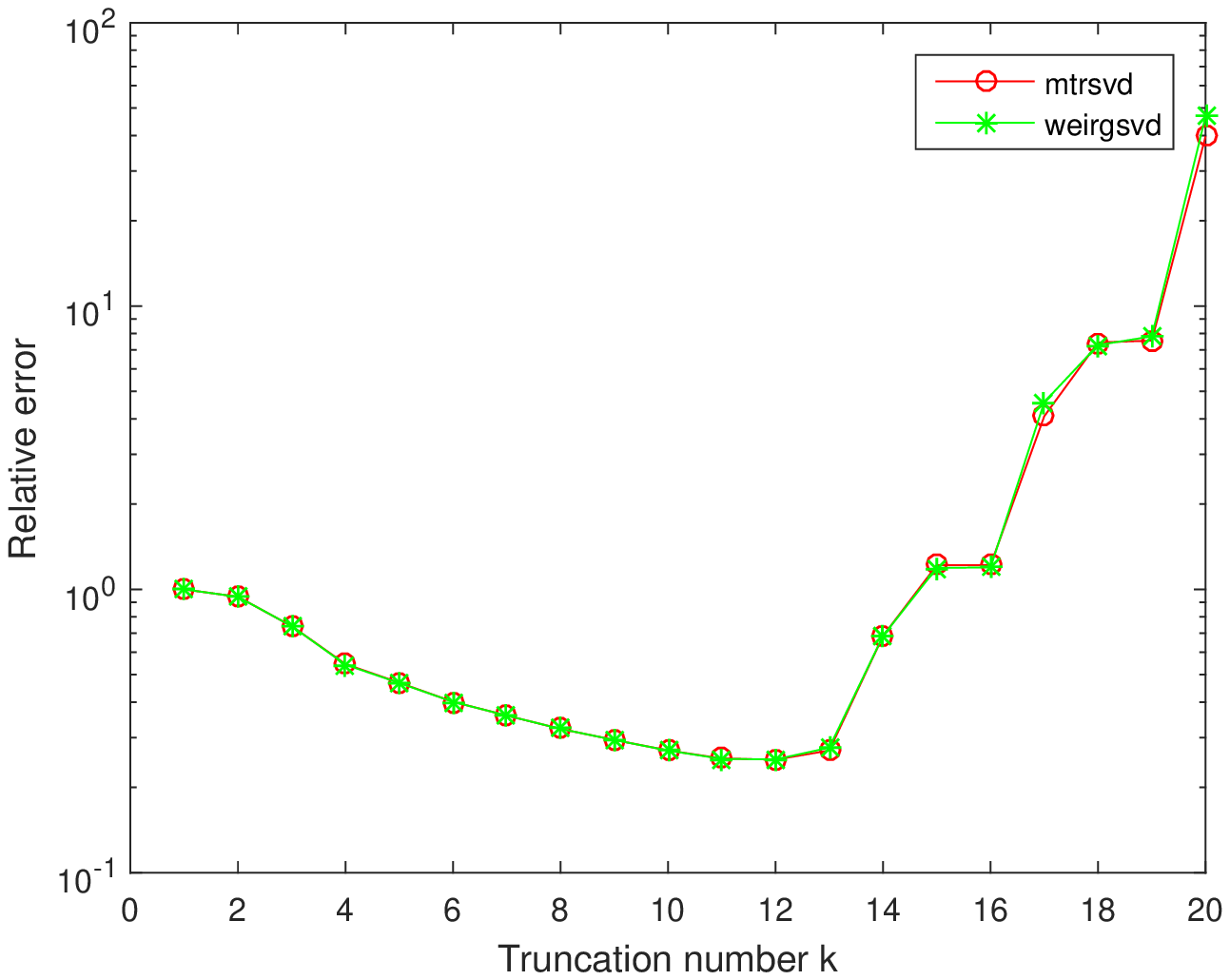}}
  \centerline{(b)}
\end{minipage}

\vfill
\begin{minipage}{0.48\linewidth}
  \centerline{\includegraphics[width=6.0cm,height=4cm]{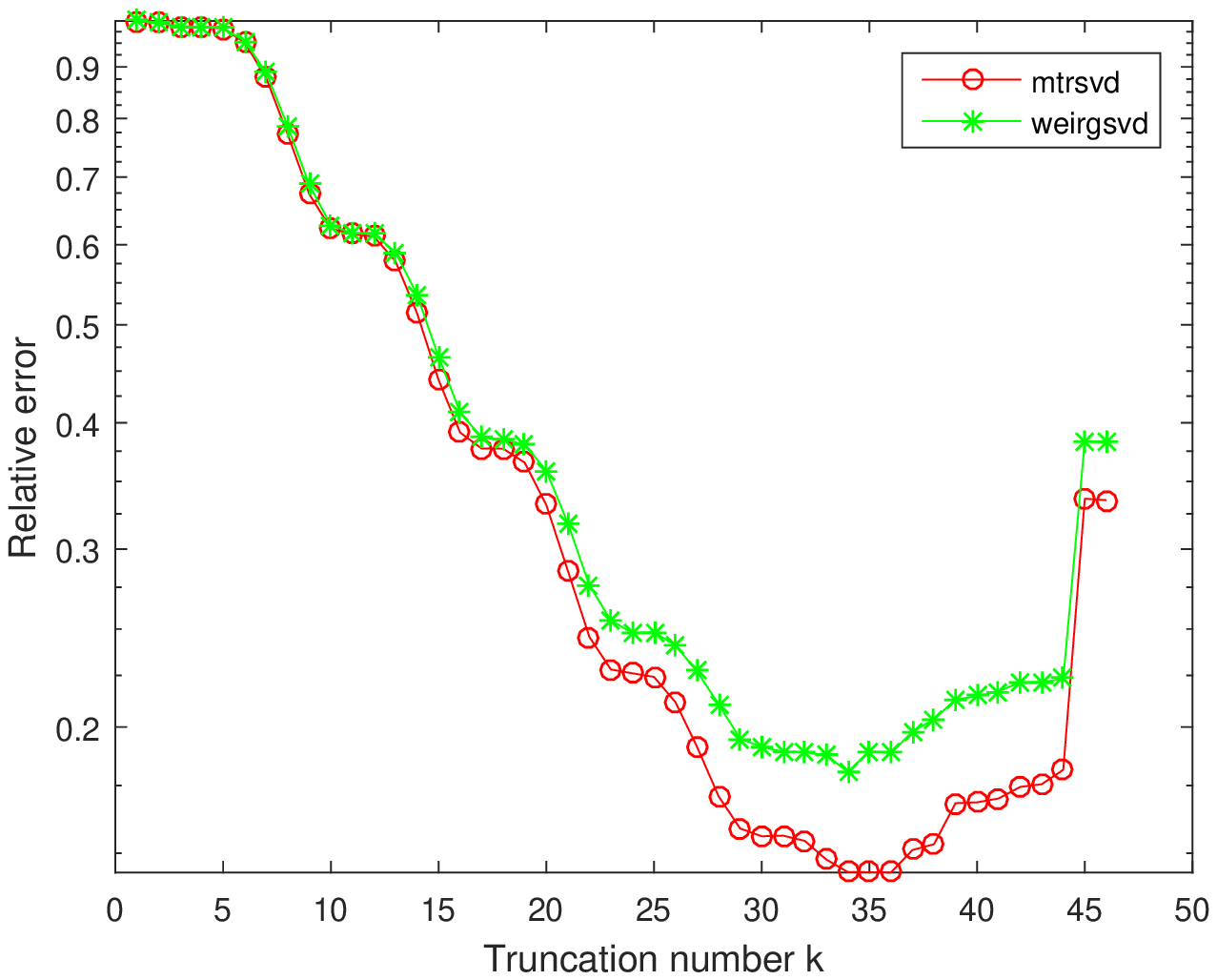}}
  \centerline{(c)}
\end{minipage}
\hfill
\begin{minipage}{0.48\linewidth}
  \centerline{\includegraphics[width=6.0cm,height=4cm]{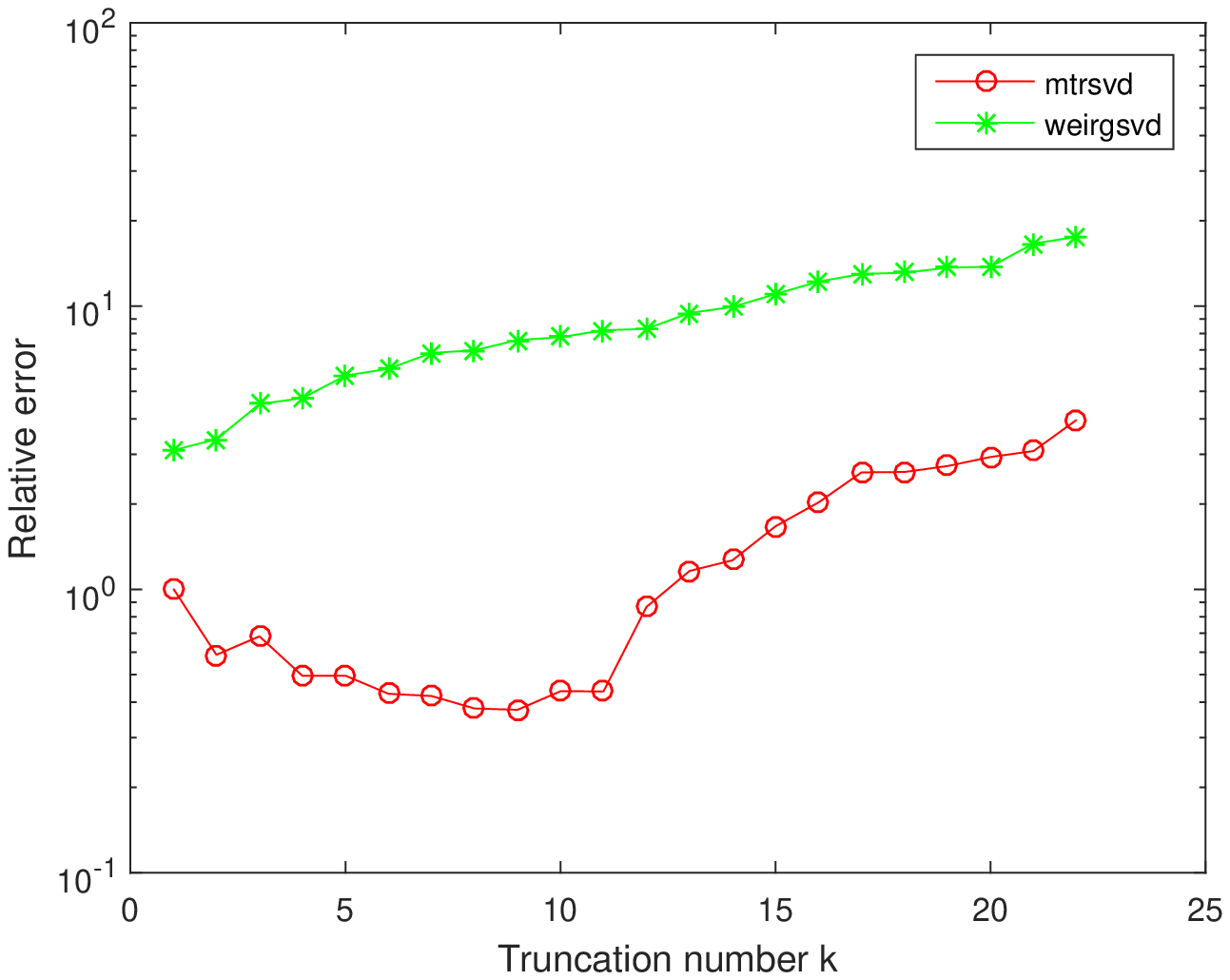}}
  \centerline{(d)}
\end{minipage}
\caption{ The relative errors of Algorithm \ref{alg5} ({\sf mtrsvd}) with $L=L_3$ and
$\varepsilon=10^{-3}$ and $m=n=10,240$: (a) {\sf shaw};
(b) {\sf gravity}; (c) {\sf heat}; (d) {\sf deriv2}.}
\label{fig6}
\end{figure}

Figure \ref{fig5} depicts the curves of convergence processes of
all the algorithms
for the four test problems {\sf shaw}, {\sf gravity}, {\sf heat} and
{\sf deriv2} with $L=L_3$, $\varepsilon=10^{-3}$ and $m=n=1,024$.
Figure \ref{fig6} does the same job for these four problems
with $L=L_3$, $\varepsilon=10^{-3}$ and $m=n=10,240$.
From the two figures, we can see that
for the severely ill-posed problem {\sf shaw} and {\sf gravity},
the relative errors obtained by
{\sf mtrsvd} are almost identical to those by {\sf tgsvd} and {\sf weirgsvd},
{\sf xiangrgsvd}.
For the moderately and mildly ill-posed problems,
{\sf mtrsvd} also behaves like {\sf tgsvd}, but
the best regularized solutions obtained by {\sf weirgsvd} and
{\sf xiangrgsvd} for the mildly ill-posed problem {\sf deriv2} have no
accuracy and their relative errors are over 300\%.
Also, we notice that for the moderately ill-posed problem {\sf heat} with
$L=L_3$, the best regularized solutions by {\sf weirgsvd} and {\sf xiangrgsvd}
are much less accurate than those by {\sf tgsvd} and {\sf mtrsvd}.

\begin{figure}
\begin{minipage}{0.48\linewidth}
  \centerline{\includegraphics[width=6.0cm,height=4cm]{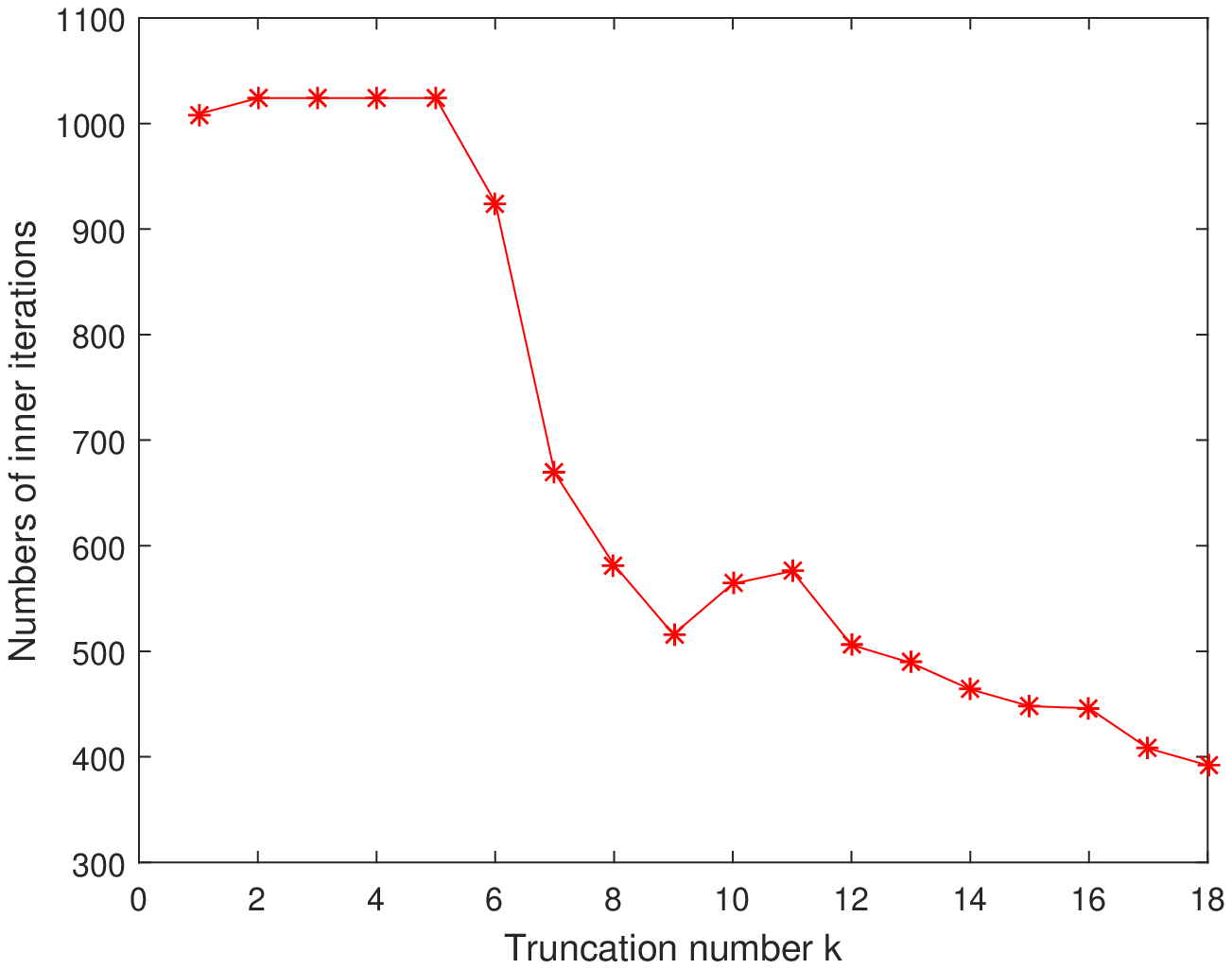}}
  \centerline{(a)}
\end{minipage}
\hfill
\begin{minipage}{0.48\linewidth}
  \centerline{\includegraphics[width=6.0cm,height=4cm]{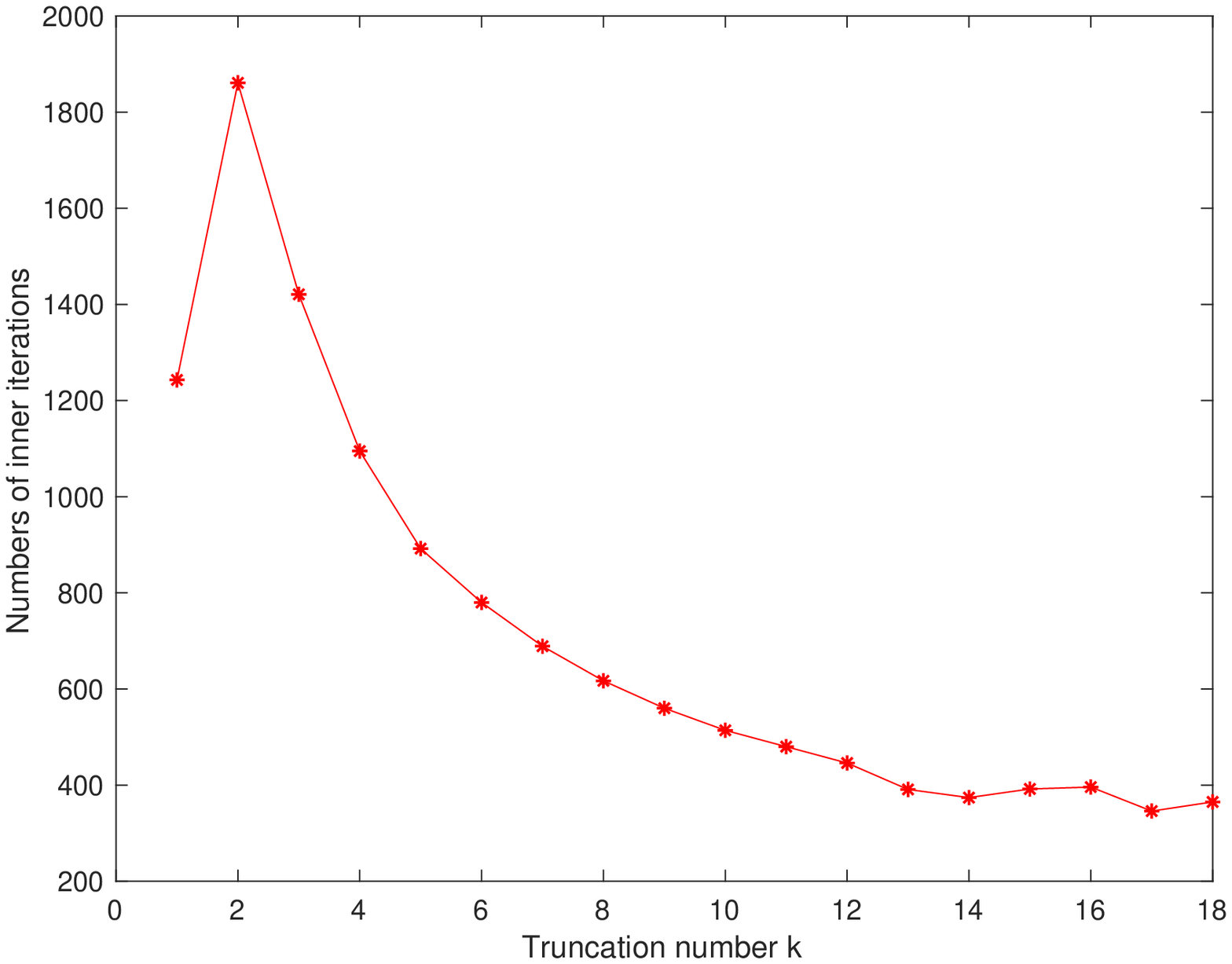}}
  \centerline{(b)}
\end{minipage}

\vfill
\begin{minipage}{0.48\linewidth}
  \centerline{\includegraphics[width=6.0cm,height=4cm]{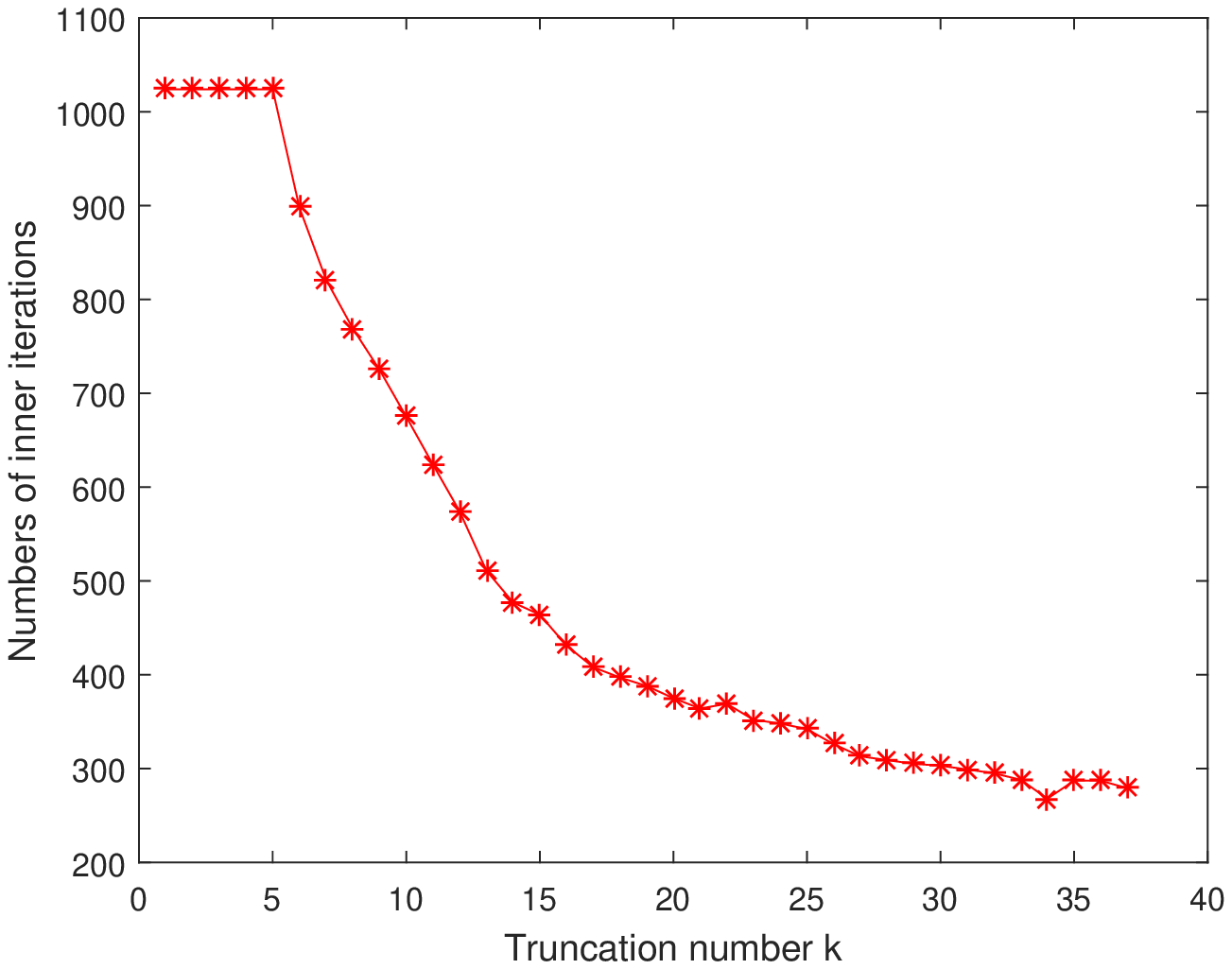}}
  \centerline{(c)}
\end{minipage}
\hfill
\begin{minipage}{0.48\linewidth}
  \centerline{\includegraphics[width=6.0cm,height=4cm]{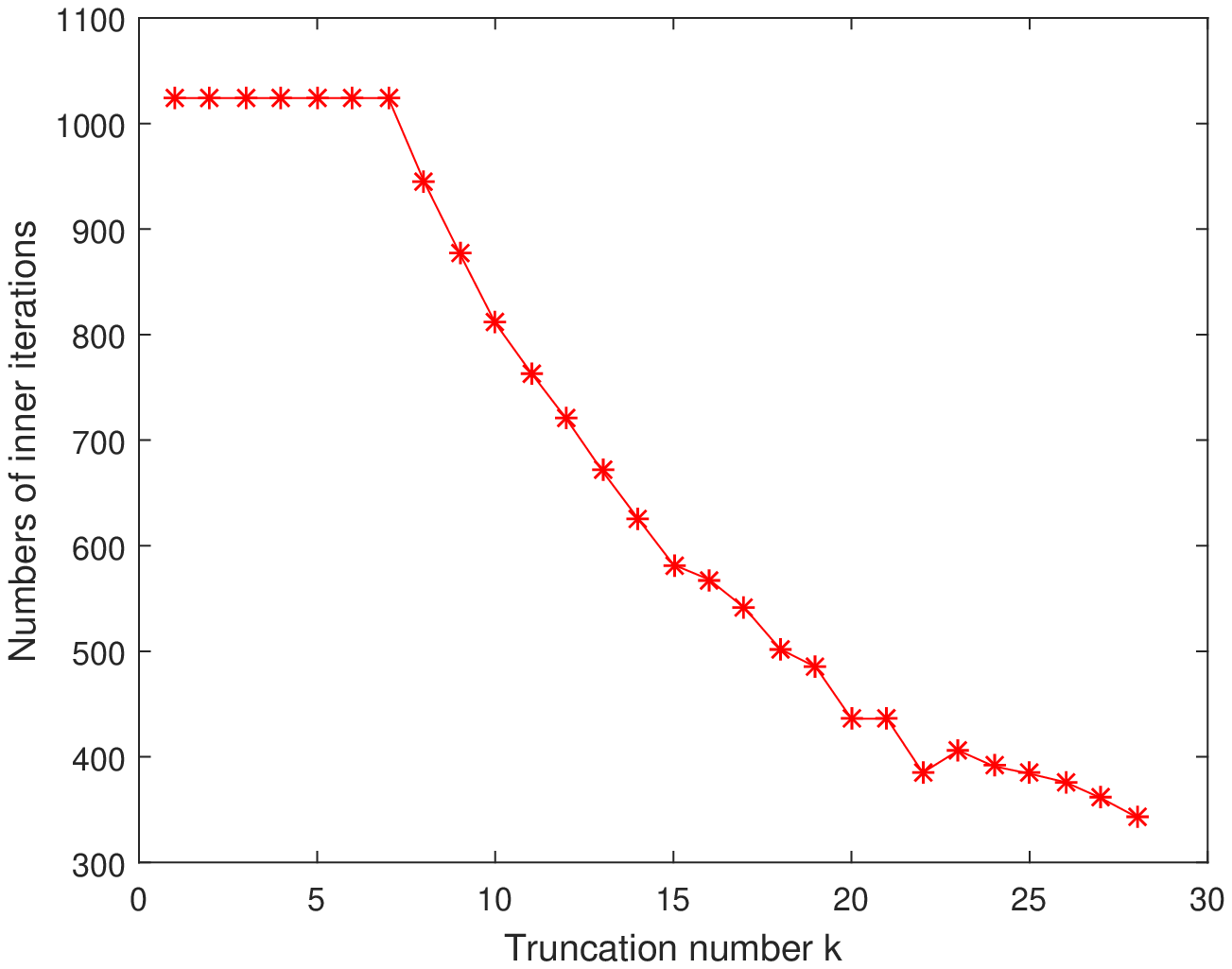}}
  \centerline{(d)}
\end{minipage}
\caption{Inner iterations versus the
truncation parameter $k$  of Algorithm~\ref{alg5} ({\sf mtrsvd})
 with $L=L_3$ and $\varepsilon=10^{-3}$ and $m=n=1,024$: (a) {\sf shaw};
 (b) {\sf gravity}; (c) {\sf heat}; (d) {\sf deriv2}.}
\label{fig7}
\end{figure}

Finally, we observe from Figure~\ref{fig7} that the number of the inner iterations
used by LSQR decrease as $k$ increases
for some chosen test problems when $L=L_3$, $\varepsilon=10^{-3}$ and $m=n=1,024$.
We see that after a few iterations, LSQR only
needs two or three hundreds iterations and even no more than one hundred
iterations to achieve the prescribed tolerance.

\section{Conclusion}

In this paper, we have proposed two MTRSVD algorithms for solving
the overdetermined
and underdetermined \eqref{eq1} with general-form regularization, respectively.
We have established a number of sharp error bounds for the approximation
accuracy of randomized approximate SVDs for three kinds of
ill-posed problems and their truncated rank-$k$ ones.
These results have improved the existing bounds substantially
and provided strong theoretical supports for the effectiveness
of randomized algorithms for solving ill-posed problems.
We have considered
the conditioning of inner least squares problems and shown
that it becomes better conditioned as the regularization
parameter $k$ increases. As a consequence, LSQR generally
converge faster with $k$ and uses fewer iterations to
achieve the prescribed tolerance, which has been confirmed numerically.
In the meantime, we have given a detailed analysis on the stopping tolerance
of LSQR for inner least squares problems and shown how to choose it
in order to guarantee that the computed regularized solutions have
the same accuracy as the ones when the problems are solved exactly.
Numerical experiments have confirmed our theory.

A practical advantage of MTRSVD is its applicability to truly large scale
problems for both overdetermined and underdetermined ill-posed
problems, while TGSVD suits only for small to medium scale problems.
For the overdetermined problems, the TRGSVD algorithm
in \cite{wei16}, though theoretically good,
are practically infeasible since it required to
compute the GSVD of the large matrix pair
$\{B,L\}$ with $B=Q^TA \in \mathbb{R}^{l \times n}$ and invert
a large $n\times n$ matrix to get the right singular vector matrix;
for the underdetermined problems, the TRGSVD algorithms in \cite{wei16,zou15}
seems to lack necessary theoretical supports and may not work well.
Some of our numerical experiments have confirmed this deficiency.

Numerical experiments have demonstrated that our MTRSVD algorithms
can compute regularized solutions with
very similar accuracy to those by the standard TGSVD algorithm and they
are at least as effective as
the TRGSVD algorithms in \cite{wei16,zou15} for solving
both overdetermined and underdetermined \eqref{eq1}.


\section*{Acknowledgements} We thank the two referees for
their comments and suggestions that further improved the presentation of our paper.

\end{document}